\documentclass[11pt,reqno]{amsart}
\usepackage{amsfonts,amsmath,amssymb}
\newtheorem{thm}{Theorem}[section]
\newtheorem{proposition}[thm]{Proposition}
\newtheorem{corollary}[thm]{Corollary}
\newtheorem{lemma}[thm]{Lemma}

\newtheorem{remark}[thm]{Remark}
\usepackage{graphicx}
\usepackage{epstopdf}

\title[Euler characters and super Jacobi polynomials] {Euler characters and super Jacobi polynomials}
\author{A.N. Sergeev}\address{Department of Mathematical Sciences,
Loughborough University, Loughborough LE11 3TU, UK}
\email{A.N.Sergeev@lboro.ac.uk}

\author{A.P. Veselov}
\address{Department of Mathematical Sciences,
Loughborough University, Loughborough LE11 3TU, UK and Moscow State University, Moscow, 119899, Russia}
\email{A.P.Veselov@lboro.ac.uk}

\begin{document}
\begin{abstract} We prove that Euler supercharacters for orthosymplectic Lie superalgebras  can be obtained as a certain specialization of super Jacobi polynomials. A new version of Weyl type formula  for super Schur functions and  specialized super Jacobi  polynomials play a key role in the proof. 
\end{abstract}

\maketitle

\tableofcontents 
\section{Introduction}

The main purpose of this paper is to develop further the link between the theory of the deformed Calogero-Moser systems and representation theory of Lie superalgebras \cite{SV,SV3}. Recall that in the $BC(m,n)$ case the deformed Calogero-Moser systems depend on 3 parameters $k,p,q.$ For generic values of these parameters they have polynomial eigenfunctions $\mathcal SJ_{\lambda}(u,v; k,p,q)$ called super Jacobi polynomials \cite{SV2}. The special case $k=p=-1, q=0$ 
corresponds to the orthosymplectic Lie superalgebra $\mathfrak{osp}(2m+1,2n)$. 

It turns out that this case is singular in the sense that the corresponding limit  does not always exist. However if we consider  first the limit when $k \rightarrow -1$ with generic $p,q$ and then let $(p,q) \rightarrow (-1,0)$ then this limit does always exist and gives what we call {\it specialized super Jacobi polynomials} $\mathcal SJ_{\lambda}(u,v; -1,-1,0).$  A natural question is what do they correspond to in the representation theory of orthosymplectic Lie superalgebras. We show that the answer is given by the so-called {\it Euler characters} studied by Penkov and Serganova \cite{PS, Ser}. 

There is a classical construction due to Borel, Weil and Bott of the irreducible representations of the complex semisimple Lie groups $G$ in terms of the cohomology of the holomorphic line bundles over the corresponding flag  varieties $G/B$ (see e.g. \cite{FH}, section 23.3). Such line bundles $L_{\lambda}$ are determined by the weight $\lambda \in \mathfrak h^{*},$ where $\mathfrak h$ is Cartan subalgebra of Lie algebra of $G$ and $B$ is Borel subgroup of $G.$ The cohomology groups $H^i(G/B, L_{\lambda})$ are finite-dimensional and have natural actions of $G$ on them. By Kodaira vanishing theorem  all of them are zero except one depending on the Weyl chamber the weight $\lambda$ belongs to (see details in \cite{FH,Bott}). In particular, for a dominant weight $\lambda$ the space of sections $H^0(G/B, L_{\lambda})$ gives the irreducible representation with highest weight $\lambda.$

In the Lie supergroup case in general there is no vanishing property, so this construction does not work \cite{PS}. The idea is to consider the virtual representation given by the Euler characteristic $$\mathcal E_{\lambda} = \sum_i (-1)^i H^i(G/B, \mathcal O_{\lambda})$$
for certain sheaf cohomology  groups (see \cite{Ser}). For the generic (typical) highest weights $\lambda$ this leads to the Kac character formula \cite{Kac}.

One can generalise this construction for any parabolic subgroup $P$ of $G$ in a natural way. In the orthosymplectic case $G=OSP(2m+1,2n)$ there is a natural choice of $P$ with the reductive part $GL(m,n).$ The corresponding supercharacter $E_{\lambda}$ (called {\it Euler supercharacter})  can be given by the general explicit formula due to Serganova \cite{Ser}. Our main result is that $E_{\lambda}$ up to a constant factor coincides with the specialized super Jacobi polynomials $\mathcal SJ_{\lambda}(u,v; -1,-1,0).$ A similar result holds for the Lie superalgebra $\mathfrak{osp}(2m,2n)$ and super Jacobi polynomials $\mathcal SJ_{\lambda}(u,v; -1,0,0).$

The proof is based on a new formula for super Schur polynomials and the super version \cite{SV2} of Okounkov's formula for Jacobi polynomials \cite{OO}. We prove also the Pieri and Jacobi--Trudy formulas for the corresponding specialized super Jacobi polynomials and Euler supercharacters.

It turns out that we can simplify the relations with super Jacobi polynomials if we choose a different  Borel subalgebra and the corresponding parabolic subalgebras following recent work  by Gruson and Serganova \cite{GS}. In that case the Euler supercharaters coincide with specialised super Jacobi polynomials without non-trivial factor (see the last section for the details). 

This shows that the super Jacobi polynomials can be considered as a natural deformation of the Euler supercharacters and gives one more evidence of a close relationship between quantum integrable systems and representation theory.

\section{Weyl type formulas for super Schur polynomials}

We start with the new formula for super Schur polynomials, which will play an important role in this work.

Let $H(m,n)$ be the set of partitions with $\lambda_{m+1} \leq n,$ which means that the corresponding Young diagram belongs to the {\it fat $(m,n)$-hook}.
Let $\lambda$ be such a partition and let $d=m-n$ be the {\it superdimension}. Introduce the following quantities 
\begin{equation}
\label{ilambda}
i(\lambda)=\max\{i\mid\lambda_{i}+d-i\ge 0,\:\; 1\le i\le m\},
\end{equation}
\begin{equation}
\label{jlambda}
j(\lambda)=\max\{j\mid\lambda^{\prime}_{j}-d-j\ge 0,\:\; 1\le j\le n\}.
\end{equation}
If all $\lambda_{i}+d-i < 0$ then by definition $i(\lambda)=0$ (and similarly for $j(\lambda)$).
It is easy to verify that in all cases $m-i(\lambda)= n-j(\lambda)$, so $i(\lambda)-j(\lambda)=d.$

We should mention that  Moens and van der Jeugt introduced a similar quantity $k(\lambda),$ which they called $(m,n)$-index of $\lambda$ (see Definition 2.2 in \cite{MV}). They were motivated by Kac-Wakimoto formula \cite{KW}. It is related to our $i(\lambda)$ by  $i(\lambda) = k(\lambda)-1.$ Moens and van der Jeugt used this quantity to write
down a new determinantal formula for super Schur polynomials different  from Sergeev-Pragacz formula \cite{PT}. 

Our formula (\ref{Weyl}) below is another new formula of Weyl type, which generalizes Sergeev-Pragacz  formula. Let us denote  by $\pi_{\lambda}$ the set of pairs $(i,j)$ such  that  $i\le i(\lambda)$ or
$j\le j(\lambda)$  and fix a partition $\nu$  such that 
$$
\lambda\cap\Pi_{m,n}\subseteq\nu\subseteq \pi_{\lambda},
$$
where $\Pi_{m,n}$ is the rectangle of the size $m\times n$. 
When $\nu = \pi_{\lambda}$ the formula (\ref{Weyl}) coincides with Serganova's formula (\ref{Vera}) for a special choice of parabolic subalgebra depending on the weight (although that was not the way we came to this).

Introduce the following quantities by
\begin{equation}
\label{frob1}
l_{i}=\lambda_{i}+m-\nu_{i}-i,\: 1\le i\le i(\lambda),\quad l_{i}=m-i, \: i(\lambda)< i\le m,
\end{equation}
\begin{equation}
\label{frob2}
k_{j}=\lambda'_{j}+n-\nu^{\prime}_{j}-j,\: 1\le j\le j(\lambda),\quad k_{j}=n-j, \: j(\lambda)< j\le n.
\end{equation}

Now we can formulate the main result of this section. Recall that {\it super Schur polynomial} is the supercharacter of the polynomial representation $M=V^{\lambda}$ of $\mathfrak{gl}(m,n)$ determined by a Young diagram $\lambda$ from the fat hook $H(m,n)$. It can be given by the following {\it Jacobi--Trudy formula} (see \cite{S}):
\begin{equation}
\label{jt}
SP_{\lambda}(x,y)=\left|\begin{array}{cccc}
   h_{\lambda_{1}}&  h_{\lambda_{1}+1}& \ldots & h_{\lambda_{1}+l-1}\\
  h_{\lambda_{2}-1}&   h_{\lambda_{2}}& \ldots &  h_{\lambda_{2}+l-2}\\
\vdots&\vdots&\ddots&\vdots\\
  h_{\lambda_{l}-l+1}& h_{\lambda_{l}-l+2} & \ldots &
  h_{\lambda_{l}}\\
 \end{array}\right|,
\end{equation} 
where  $l=l(\lambda)$ is the number of non-zero parts in the partition $\lambda$, $h_k = h_k(x,y)$
are determined by 
\begin{equation}\label{prod}
\frac{\prod_{j=1}^n(1-ty_{j})}{\prod_{i=1}^m(1-tx_{i})}=\sum_{a=0}^{\infty}h_{a}(x,y)t^a
\end{equation}
and $h_{a}=0$ if $a<0$. One can check that the highest coefficient  of $SP_{\lambda}(x,y)$ is equal to
$(-1)^b,$ where $b=\sum_{j>m}\lambda_{j}.$ This explains the appearance of this sign below (see also \cite{SV1}).

\begin{thm} \label{SS} The super  Schur polynomial $SP_{\lambda}(x_{1},\dots,x_{m},y_{1},\dots,y_{n})$ can be expressed by the following Weyl type formula for any choice of partition $\nu$ such that $\lambda\cap\Pi_{m,n}\subseteq\nu\subseteq \pi_{\lambda}$:
\begin{equation}\label{Weyl}
\begin{split}
SP_{\lambda}&(x_{1},\dots,x_{m},y_{1},\dots,y_{n})=\\
&(-1)^{b}\sum_{w\in S_{m}\times S_{n}}w\left[\prod_{(i,j)\in\nu}(x_{i}-y_{j})\frac{x_{1}^{l_{1}}\dots x_{m}^{l_{m}}y_{1}^{k_{1}}\dots y_{n}^{k_{n}}}{\Delta(x)\Delta(y)}\right],
\end{split}
\end{equation}
where $\Delta(x) = \prod_{i<j}^m(x_i-x_j),\, \Delta(y)=\prod_{i<j}^n(y_i-y_j), \, b=\sum_{j>m}\lambda_{j}.$
\end{thm}
\begin{proof}
For any function $f(x,y)$ define the following alternation operations
$$
\{f(x,y)\}=\sum_{w\in S_{m}\times S_{n}}\varepsilon(w)w(f(x,y))
$$
and
$$
\{f(x,y)\}_{x}=\sum_{w\in S_{m}}\varepsilon(w)w(f(x,y)),
$$
where $S_m$ and $S_n$ permute $x_i$ and $y_j$ respectively.
Introduce also the notations
$$x^{\rho_m}=x_{1}^{m-1}x_{2}^{m-2}\dots x_{m}^0, \quad y^{\rho_n}=y_{1}^{n-1}y_{2}^{n-2}\dots y_{n}^0.$$
First let us prove the following equality 
\begin{equation}\label{h}
\{h_{a}(x,y)x^{\rho_m}y^{\rho_n}\}=\left\{\prod_{j=1}^n(x_{1}-y_{j})x_{1}^{a+d-1}x_{2}^{m-2}\dots x_{m}^0y_{1}^{n-1}\dots y_{n}^0\right\},
\end{equation}
where $a$ is an integer such that $a+d-1\ge 0.$ 
Indeed, we have from the usual Weyl formula for $a\ge0$
$$
\{h_{a}(x)x_{1}^{m-1}x_{2}^{m-2}\dots x_{m}^0\}_{x}=\{x_{1}^{a+m-1}x_{2}^{m-2}\dots x_{m}^0\}_{x}
$$ 
This is true also for all $a\ge 1-m$ because the left hand side for negative $a$ is zero by definition.
From (\ref{prod}) we have 
$$
h_{k}(x,y)=\sum_{j=0}^{n}(-1)^{j}h_{k-j}(x)e_{j}(y),
$$
where $h_k$ and $e_j$ are complete symmetric and elementary symmetric polynomials respectively.
Note now that if $a+d-1\ge0$ and $0\le j\le n$ we have $a-j\ge 1-m$. Therefore in that case
$$
\{h_{a}(x,y)x_{1}^{m-1}\dots x_{m}^0\}_{x}=\sum_{j=0}^{n}\{h_{a-j}(x)x_{1}^{m-1}\dots x_{m}^0\}_{x}(-1)^{j}e_{j}(y)=
$$
$$
\sum_{j=0}^{n}\{x_{1}^{a-j}x^{\rho_m}\}_{x}(-1)^{j}e_{j}(y)=
\{\sum_{j=0}^{n}x_{1}^{n-j}(-1)^{j}e_{j}(y)x_{1}^{a+d-1}x_2^{m-2}\dots x_{m}^0\}_{x}=
$$
$$
\{\prod_{j=1}^n(x_{1}-y_{j})x_{1}^{a+d-1}x_2^{m-2}\dots x_{m}^0\}_{x}.
$$
This implies the formula (\ref{h}).

We prove now the theorem by induction. For this we need the following 
\begin{lemma}
\label{red} If $\lambda_{1}+d-1\ge0$ we have the following equality
\begin{equation}
\label{sp}
\left\{SP_{\lambda}(x,y)x^{\rho_m}y^{\rho_n}\right\}=
\left\{\prod_{j=1}^n(x_{1}-y_{j})x_{1}^{\lambda_{1}+d-1}SP_{\hat\lambda}(\hat x, y)x_{2}^{m-2}\dots x_{m}^0y^{\rho_n}\right\}
\end{equation}
where $\hat x= (x_2, \dots, x_m)$, $\hat\lambda=(\lambda_{2},\lambda_{3},\dots).$ 
\end{lemma}
\begin{proof}
To prove it we use the Jacobi-Trudy formula  (\ref{jt}).
Let $\lambda_{1}+d-1\ge0,$ then we have
$$
\left\{SP_{\lambda}(x_{1},\dots,x_{m},y_{1},\dots,y_{n})x_{1}^{m-1}\dots x_{m}^0y_{1}^{n-1}\dots y_{n}^0\right\}=
$$
$$
\left\{\left|\begin{array}{cccc}
   h_{\lambda_{1}}&  h_{\lambda_{1}+1}& \ldots & h_{\lambda_{1}+l-1}\\
  h_{\lambda_{2}-1}&   h_{\lambda_{2}}& \ldots &  h_{\lambda_{2}+l-2}\\
\vdots&\vdots&\ddots&\vdots\\
  h_{\lambda_{l}-l+1}& h_{\lambda_{l}-l+2} & \ldots &
  h_{\lambda_{l}}\\
 \end{array}\right|x_{1}^{m-1}\dots x_{m}^0y_{1}^{n-1}\dots y_{n}^0\right\}=
$$
$$
\left|\begin{array}{cccc}
   \left\{h_{\lambda_{1}}x_{1}^{m-1}\dots  y_{n}^0\right\}&  \left\{h_{\lambda_{1}+1}x_{1}^{m-1}\dots  y_{n}^0\right\}& \ldots & \left\{h_{\lambda_{1}+l-1}x_{1}^{m-1}\dots  y_{n}^0\right\}\\
  h_{\lambda_{2}-1}&   h_{\lambda_{2}}& \ldots &  h_{\lambda_{2}+l-2}\\
\vdots&\vdots&\ddots&\vdots\\
  h_{\lambda_{l}-l+1}& h_{\lambda_{l}-l+2} & \ldots &
  h_{\lambda_{l}}\\
 \end{array}\right|=
$$
\medskip
$$
\left\{\prod_{j=1}^n(x_{1}-y_{j})\left|\begin{array}{cccc}
   x_{1}^{\lambda_{1}+d-1}\dots  y_{n}^0&  x_{1}^{\lambda_{1}+d}\dots  y_{n}^0& \ldots & x_{1}^{\lambda_{1}+d+l-2}\dots  y_{n}^0\\
  h_{\lambda_{2}-1}&   h_{\lambda_{2}}& \ldots &  h_{\lambda_{2}+l-2}\\
\vdots&\vdots&\ddots&\vdots\\
  h_{\lambda_{l}-l+1}& h_{\lambda_{l}-l+2} & \ldots &
  h_{\lambda_{l}}\\
 \end{array}\right|\right\}.
$$
Now multiplying every column except  the last one by $x_{1}$ and subtracting it from the next column taking into account the equality 
$$
h_{a-1}(x,y)-x_{1}h_{a}(x,y)=h_{a}(\hat x,y)
$$
we get
$$
\left\{\prod_{j=1}^n(x_{1}-y_{j}) x_{1}^{\lambda_{1}+d-1}\left|\begin{array}{cccc}
  \hat h_{\lambda_{2}}&   \hat h_{\lambda_{2}+1}& \ldots &  \hat h_{\lambda_{2}+l-2}\\
\vdots&\vdots&\ddots&\vdots\\
  \hat h_{\lambda_{l}-l+2}&\hat  h_{\lambda_{l}-l+3} & \ldots &
 \hat h_{\lambda_{l}}\\
 \end{array}\right|\right\},
$$ 
where $\hat h_k = h_k(\hat x,y)$ and $\hat x =(x_2, \dots, x_m)$ as before.
This implies the formula (\ref{sp}).
\end{proof}

\begin{lemma}\label{stab} Let $\lambda_{1} < n$, then we have the following equality
$$
\left\{SP_{\lambda}(x,y_{1},\dots,y_{n})x^{\rho_m}y^{\rho_n}\right\}=\left\{SP_{\lambda}(x,y_{1},\dots,y_{n-1})x^{\rho_m}y^{\rho_n}\right\}.
$$
\end{lemma}
\begin{proof} We have (see formulas (5.9) in \cite{Mac} and (25) in \cite{SV1})
$$
SP_{\lambda}(x,y)=\sum_{\mu \subseteq\lambda}(-1)^{|\mu|}S_{\lambda/\mu}(x)S_{\mu^{\prime}}(y)
$$
So, if $\lambda_{1} < n$ then $\mu_{n}^{\prime}=0$ and lemma follows.
\end{proof}
Now we finish the proof  by induction in $m+n.$
Note first that the case $n=0$ follows from the classical Weyl's formula, while in the case $m=0$ we have also the sign $(-1)^{\lambda_{m+1}+\lambda_{m+2}+\dots}$ determined as in the formula (25) in \cite{SV1}. Let us assume now that $m>0.$
Rewrite formula (\ref{Weyl}) in the form
\begin{equation}\label{Weyl'}
\left\{SP_{\lambda}(x,y)x^{\rho_m}y^{\rho_n}\right\}=
(-1)^{b}\left\{\prod_{(i,j)\in\nu}(x_{i}-y_{j}) x_{1}^{l_{1}}\dots x_{m}^{l_{m}}y_{1}^{k_{1}}\dots y_{n}^{k_{n}}\right\}.
\end{equation}
If $\nu_1 < n$ (and hence $\lambda_{1} < n$) by induction assumption we have
\begin{equation*}
\begin{split}
&\left\{SP_{\lambda}(x,y_{1},\dots,y_{n-1})x^{\rho_m}y^{\rho_{n-1}}\right\}=\\
&(-1)^{b}\left\{\prod_{(i,j)\in\nu}(x_{i}-y_{j}) x_{1}^{l_{1}}\dots x_{m}^{l_{m}}y_{1}^{k_{1}-1}\dots y_{n-1}^{k_{n-1}-1}\right\}.
\end{split}
\end{equation*}
Multiplying both sides by the product $y_1\dots y_n$ and using lemma \ref{stab}, we prove the theorem in this case.

If $\nu_1=n$ then the box $(1,n)$ belongs to $\pi_{\lambda}$ and therefore by definition either $i(\lambda)>0$ or $j(\lambda)=n.$ If $i(\lambda)=0$ then $j(\lambda)=n-m<n$ since $m>0.$ Thus, $i(\lambda)\ge 1,$
which implies that $\lambda_1+d-1 \ge 0.$ This means that we can apply lemma \ref{red}.
By induction we have
\begin{equation*}
\begin{split}
&\left\{SP_{\lambda}(x_{2},\dots,x_{m})x_2^{m-2}\dots x_m^0 y^{\rho_n}\right\}=\\
&(-1)^{b}\left\{\prod_{(i,j)\in\hat \nu}(x_{i}-y_{j}) x_{2}^{l_{2}}\dots x_{m}^{l_{m}}y_{1}^{k_{1}}\dots y_{n}^{k_{n}}\right\},
\end{split}
\end{equation*}
where $\hat\nu$ is the partition $\nu$ without the first row. Now the theorem follows from lemma \ref{red}.
 \end{proof}

\section{Euler supercharacters for Lie superalgebra $\mathfrak{osp}(2m+1,2n)$}

In this section we consider the case of orthosymplectic Lie superalgebra $\mathfrak{osp}(2m+1,2n),$
the case of $\mathfrak{osp}(2m,2n)$ is considered in section 7.

Recall the description of the root system of Lie superalgebra $\mathfrak{g}=\mathfrak{osp}(2m+1,2n)$.
We have $\mathfrak{g}= \mathfrak{g}_0 \oplus \mathfrak{g}_1$, where $\mathfrak{g}_{0}=so(2m+1)\oplus sp(2n)$ and
$\mathfrak{g}_{1}=V_{1}\otimes V_{2}$ where $V_{1}$ and $V_{2}$ are the
identical representations of $so(2m+1)$ and $sp(2n)$ respectively. Let
$\pm\varepsilon_{1},\dots,\pm\varepsilon_{m},\, \pm\delta_1, \dots, \pm\delta_n$ be the non-zero weights of the
identical representation of  $\mathfrak{g}.$
The root system of $osp(2m+1,2n)$ consists of
$$R_{0}=\{\pm\varepsilon_{i}\pm\varepsilon_{j},\:\pm\varepsilon_{i}, \, \pm\delta_{p}\pm\delta_{q},\:\pm2\delta_{p}, \,
\: i\ne j,\,\: 1\le i,j\le m\ , p\ne q,\, 1\le p,q\le n\},
$$
$$
R_{1}=\{\pm\varepsilon_{i}\pm\delta_{p},\:\pm\delta_{p}\},\quad
R_{iso}=\{\pm\varepsilon_{i}\pm\delta_{p}\},
$$
where $R_{0}, R_{1}$ and $R_{iso}$ are even, odd and isotropic parts respectively.
The invariant bilinear form is given by
$$
(\varepsilon_{i},\varepsilon_{i})=1,\: (\varepsilon_{i},\varepsilon_{j})=0,\: i\ne j,\:
(\delta_{p},\delta_{p})=-1,\: (\delta_{p},\delta_{q})=0,\: p\ne q,\: (\varepsilon_{i},\delta_{p})=0.
$$
The Weyl group $W_0=\left(S_{m}\ltimes\mathbb{Z}_{2}^m\right)\times\left( S_{n}\ltimes\mathbb{Z}_{2}^n\right)$ acts on the weights by separately  permuting $\varepsilon_{i},\; j=1,\dots,m$ and
$\delta_{p},\; p=1,\dots,n$ and changing their signs.
A distinguished system of simple roots can be chosen as
$$
B=\{\delta_{1}-\delta_{2},\dots,\delta_{n-1}-\delta_{n},\delta_{n}-\varepsilon_{1},\varepsilon_{1}-\varepsilon_{2},\dots,\varepsilon_{m-1}-\varepsilon_{m}, \varepsilon_{m}\}.
$$
Introduce the variables $x_{i}=e^{\varepsilon_{i}},\:x^{1/2}_{i}=e^{\varepsilon_{i}/2},\, u_{i}=x_{i}+x_{i}^{-1},\, i=1,\dots,m$ and $ y_{p}=e^{\delta_{p}},\:  v_{p}=y_{p}+y_{p}^{-1},\: p=1,\dots, n$.

For any parabolic subalgebra $\mathfrak p \subset \mathfrak g$ and any finite-dimensional representation $M$ of $\mathfrak p$ by a superversion of Borel-Weil-Bott construction one can define the corresponding Euler supercharacter $E^{\mathfrak p}(M)$. 
According to the general formula due to Serganova \cite{Ser}
\begin{equation}
\label{Vera}
E^{\mathfrak p}(M) = \sum_{w \in W_0} \, w\left(\frac{D \, e^{\rho} sch M}{\prod_{\alpha \in R_{\mathfrak p} \cap R_1^+} (1-e^{-\alpha})}\right)
\end{equation}
with
$$D= \frac{\prod_{\alpha \in R_1^+}(e^{\alpha/2} - e^{-\alpha/2})}{\prod_{\alpha \in R_0^+}(e^{\alpha/2} - e^{-\alpha/2})}.$$
Here $\rho$ is the half-sum of the even positive roots minus the half-sum of odd positive roots, $R_{\mathfrak p}$ is the set of roots $\alpha$ such that $\mathfrak g_{\pm\alpha} \subset \mathfrak p$
(see formula (3.1) in \cite{Ser}). Note that we use here the supercharacter rather than the character used by Serganova. This leads simply to the change of signs in some places.

Consider now the parabolic subalgebra $\mathfrak p$ with 
$$R_{\mathfrak p} = \{\varepsilon_{i}-\varepsilon_{j}, \, \delta_{p}-\delta_{q},\,\, \pm(\varepsilon_{i}- \delta_{p})\}.$$
The algebra $\mathfrak p$ is isomorphic to the sum of the Lie superalgebra $\mathfrak{gl}(m,n)$ and some nilpotent Lie superalgebra. In that case Serganova's formula (\ref{Vera}) has the form
\begin{equation}
\label{Vera1}
E(M) = \sum_{w \in W_0} w \left( \frac{\displaystyle\prod_{(i,j)\in \Pi_{m,n}}(1-x_i^{-1}y_j^{-1}) x_1^{m-\frac{1}{2}}\dots x_m^{\frac{1}{2}} y_1^{n-\frac{1}{2}}\dots y_n^{\frac{1}{2}} sch M}{\Delta(u)\Delta(v)\prod_{i=1}^m (x_i^{\frac{1}{2}}-x_i^{-\frac{1}{2}}) \prod_{j=1}^n (y_j^{\frac{1}{2}}+y_j^{-\frac{1}{2}})}\right),
\end{equation}
where $\Pi_{m,n}$ is the rectangle of the size $m\times n$, $\Delta(u) = \prod_{i<j}^m(u_i-u_j),\, \Delta(v)=\prod_{i<j}^n(v_i-v_j)$ and $u_{i}=x_{i}+x_{i}^{-1}, v_{j}=y_{j}+y_{j}^{-1}$ as before.

The following proposition explains the appearance of the powers of 2 in the later considerations. As far as we know the appearance of a power of 2 coefficient was first noticed by Cheng and Wang \cite{CW}.
It should have a geometrical explanation but we give here a direct algebraic proof. 

\begin{proposition} \label{Triv} For the trivial even representation $M$ we have
\begin{equation}
\label{triv}
E(M) = 2^{\min(m,n)}.
\end{equation}
\end{proposition}

\begin{proof} It is well-known (see e.g. \cite{Mac}, Ch. 1, formula (4.3')) that
$$\prod_{i,j}(1-x_i^{-1}y_j^{-1})= \sum (-1)^{|\lambda|} S_{\lambda}(x_1^{-1},\dots, x_m^{-1})S_{\lambda'}(y_1^{-1},\dots, y_n^{-1}),$$
where $S_{\lambda}$ is the Schur polynomial  and the sum is over all the Young diagrams $\lambda$, which are contained in the $(m\times n)$ rectangle, $\lambda'$ is the diagram transposed to $\lambda.$
Using the Weyl formula for Schur polynomials we can replace in the formula (\ref{Vera1}) the product
$\prod_{i,j}(1-x_i^{-1}y_j^{-1})$ by the sum $$\sum (-1)^{|\lambda|} x_1^{-\lambda_m+m-1/2}\dots x_m^{-\lambda_1+1/2} y_1^{-\lambda'_n+n-1/2}\dots y_n^{-\lambda'_1+1/2}.$$ One can check that the only non-vanishing terms correspond to the symmetric Young diagrams $\lambda = \lambda'.$
Now the proposition follows from the fact that the number of the symmetric diagrams contained in the $(m\times n)$ rectangle is equal to $2^{\min(m,n)},$ which can be easily proved by induction.
\end{proof}

Consider now the polynomial representations $M=V^{\lambda}$ of $\mathfrak{gl}(m,n)$ determined by a Young diagram $\lambda \in H(m,n)$ with the supercharacter given by the super Schur polynomial:
$$sch \, V_{\lambda} = SP_{\lambda}(x_1, \dots,x_m, y_1,\dots,y_n).$$
The next step is to rewrite the formula for the Euler supercharacter
\begin{equation}
\label{el}
E_{\lambda} = E(V^{\lambda}),
\end{equation}
 using the formula (\ref{Weyl}) for the super Schur polynomials in the case when $\nu = \pi_{\lambda}.$
 In this case $k_i, l_j$ are defined by 
 \begin{equation}
\label{frob11}
l_{i}=\lambda_{i}+d-i,\: 1\le i\le i(\lambda),\quad l_{i}=m-i, \: i(\lambda)< i\le m,
\end{equation}
\begin{equation}
\label{frob21}
k_{j}=\lambda'_{j}-d-j,\: 1\le j\le j(\lambda),\quad k_{j}=n-j, \: j(\lambda)< j\le n,
\end{equation}
where $d=m-n.$

Introduce the following polynomials (which are particular cases of classical Jacobi polynomials, see the next section)
\begin{equation}
\label{jaco}
\varphi_a(z) = \frac{w^{a+1/2} - w^{-a-1/2}}{w^{1/2} - w^{-1/2}}, \quad 
\psi_a(z) = \frac{w^{a+1/2} + w^{-a-1/2}}{w^{1/2} + w^{-1/2}},
\end{equation}
where $z=w+w^{-1}.$ Define also $\Pi_{\lambda}(u,v)$ as
$$\Pi_{\lambda}(u,v)= \prod_{(i,j)\in \pi_{\lambda}}(u_i-v_j).$$

\begin{thm} The Euler supercharacters can be given by the following formula
\begin{equation}\label{Euler}
\begin{split}
&E_{\lambda}(u,v)= \\
&C(\lambda)\sum_{w\in S_{m}\times S_{n}}w\left[\Pi_{\lambda}(u,v)\frac{\varphi_{l_{1}}(u_1)\dots {\varphi_{l_{m}}(u_m)\psi_{k_{1}}(v_1)\dots \psi_{k_{n}}(v_n)}}{\Delta(u)\Delta(v)}\right],
\end{split}
\end{equation}
where $C(\lambda) = (-1)^{b} 2^{m-i(\lambda)} = (-1)^{b} 2^{n-j(\lambda)},$ $b=\sum_{j>m}\lambda_{j}$
and $i(\lambda),\, j(\lambda),\, l_{i},\, k_j$ are defined by 
(\ref{ilambda}), (\ref{jlambda}),(\ref{frob11}), (\ref{frob21}). 
\end{thm}

 \begin{proof}
According to Theorem \ref{SS} for $M=V^{\lambda}$ the corresponding supercharacter $sch M=SP_{\lambda}$  can be given by (\ref{Weyl}), or, equivalently, by (\ref{Weyl'}). 
 Substituting this into Serganova's formula (\ref{Vera1}) we get
$$
E(M) = (-1)^{b}\sum_{w \in W_0} w \left( \frac{\displaystyle\prod_{(i,j)\in \Pi_{m,n}}(1-x_i^{-1}y_j^{-1}) \prod_{(i,j)\in\pi_{\lambda}}(x_{i}-y_{j}) x^{l+\frac{1}{2}} y^{k+\frac{1}{2}}}{\Delta(u)\Delta(v)\prod_{i=1}^m (x_i^{\frac{1}{2}}-x_i^{-\frac{1}{2}}) \prod_{j=1}^n (y_j^{\frac{1}{2}}+y_j^{-\frac{1}{2}})}\right),
$$
where we use the notations
$$x^{l+\frac{1}{2}}=x_{1}^{l_{1}+\frac{1}{2}}\dots x_{m}^{l_{m}+\frac{1}{2}}, \, y^{k+\frac{1}{2}}= y_{1}^{k_{1}+\frac{1}{2}}\dots y_{n}^{k_{n}+\frac{1}{2}}.$$
Since 
$$\displaystyle\prod_{(i,j)\in \pi_{\lambda}}(1-x_i^{-1}y_j^{-1}) \prod_{(i,j)\in \pi_{\lambda}}(x_{i}-y_{j})=
\prod_{(i,j)\in \pi_{\lambda}}(u_i-v_j)=\Pi_{\lambda}(u,v),$$
we have
$$
E(M) = (-1)^{b}\sum_{w \in W_0} w \left( \frac{\Pi_{\lambda}(u,v)\displaystyle\prod_{(i,j)\in \Pi_{m,n}\setminus \pi_{\lambda}} (1-x_i^{-1}y_j^{-1})  x^{l+\frac{1}{2}} y^{k+\frac{1}{2}}}{\Delta(u)\Delta(v)\prod_{i=1}^m (x_i^{\frac{1}{2}}-x_i^{-\frac{1}{2}}) \prod_{j=1}^n (y_j^{\frac{1}{2}}+y_j^{-\frac{1}{2}})}\right).
$$
Now proposition \ref{Triv} and Weyl's formula for the $BC_n$ root system (which is the root system of the Lie superalgebra $\mathfrak{osp}(1,2n)$) allow us to replace here $\prod_{(i,j)\in \Pi_{m,n}\setminus \pi_{\lambda}} (1-x_i^{-1}y_j^{-1})$ by $2^{m-i(\lambda)}= 2^{n-j(\lambda)}$ to come to
\begin{equation}
\label{bcd}
E(M) = C(\lambda)\sum_{w \in W_0} w \left( \frac{\Pi_{\lambda}(u,v) x^{l+\frac{1}{2}} y^{k+\frac{1}{2}}}{\Delta(u)\Delta(v){\displaystyle \prod_{i=1}^m (x_i^{\frac{1}{2}}-x_i^{-\frac{1}{2}}) \prod_{j=1}^n (y_j^{\frac{1}{2}}+y_j^{-\frac{1}{2}})}}\right).
\end{equation}
Summing now over the subgroup $\mathbb Z_2^m \times \mathbb Z_2^n \subset W_0$ and using (\ref{jaco}) we have the claim.
 \end{proof}

We will use this formula now to show the relation with super Jacobi polynomials.

\section{Super Jacobi polynomials for $k=-1$}

The main result of this section is the following Weyl-type formula for the super Jacobi polynomials \cite{SV2} with $k=-1.$ Let us introduce the following polynomials $f_l(z,p,q)$, which are certain normalised versions of the classical Jacobi polynomials $P_l^{\alpha,\beta}(z)$ with $\alpha=-p-q-\frac{1}{2},\, \beta =q-\frac{1}{2}$:
$$
f_l(z,p,q)=\sum_{i=0}^l C_{l,i}(z-2)^i,
$$
where $C_{l,l}=1$ and 
$$
C_{l,i}=4^{l-i}\frac{(i+1)\dots(l-1) l}{(l-i)!}\frac{(i+1-p-q-1/2)\dots (l-p-q-1/2)}{(l+i-p-2q)\dots (2l-1-p-2q)}
$$
for $i<l.$
Introduce also $$g_{k}(w,p,q)=f_k(w,-p,-1-q).$$ Note that the polynomials $\varphi(z), \psi(z)$ from the previous section are the particular cases:
$$\varphi_l(z) = f_l(z,-1,0), \quad \psi_k(z) = g_k(z,-1,0)=f_k(z,1,-1).$$ 
Later we drop the parameters for brevity, writing simply
$f_k(z), g_k(z).$

Let $\lambda \in H(m,n)$ be a partition from the fat hook  and  $\pi_{\lambda},\, l_i, k_j$ be the same as in the previous section (see formulas (\ref{frob11}),(\ref{frob21}) above).

\begin{thm} \label{main}  The super Jacobi polynomials for special value of parameter $k=-1$ can be given by the following formula
\begin{equation}\label{WeylJ}
\begin{split}
&SJ_{\lambda}(u,v,-1,p,q)=\\
&(-1)^b\sum_{w\in S_{m}\times S_{n}}w\left[\frac{\Pi_{\lambda}(u,v)} {\Delta(u)\Delta(v)}f_{l_{1}}(u_{1})\dots f_{l_{m}}(u_{m})g_{k_{1}}(v_{1})\dots g_{k_{n}}(v_{n})\right],
\end{split}
\end{equation}
where as before $b=\sum_{j>m}\lambda_j$  and $\Pi_{\lambda}(u,v)=\prod_{(i,j)\in\pi_{\lambda}}(u_{i}-v_{j}).$
\end{thm}

We prove this first in the particular case when $\lambda$ contains the  $m \times n$ rectangle, i.e.
$$\lambda_m \geq n.$$ The corresponding formula for super Jacobi polynomials can be considered as a natural analogue of Berele-Regev factorisation formula for super Schur polynomials \cite{BR}. 
For such a diagram $\lambda$ one can consider its sub-diagram $\mu,$ which is the diagram $\lambda$ without first $n$ columns. Define
\begin{equation}
\label{embedding}
w_i(\lambda)=\mu_{i},\: i=1,\dots,m, \,\, z_j(\lambda)=\lambda^{\prime}_{j},\: j=1,\dots,n.
\end{equation}
In other words, $w_i$ is the length of $i$-th row of $\mu$ and $z_j$ is the length of $j$-th column of $\lambda \in H_{m,n}.$

The super Jacobi polynomials \cite{SV2} in the special case $k=-1$ can be defined for generic $p,q$
in terms of super Schur polynomials $SP_{\lambda}$ by Okounkov's formula 
\begin{equation}\label{okoun}
SJ_{\lambda}(u, v,-1,p,q)=\sum_{\tilde\lambda\subseteq\lambda}K_{\lambda,\tilde\lambda}SP_{\tilde\lambda}(\hat u, \hat v)
\end{equation}
where
$u=(u_{1},\dots,u_{m}), \, v=(v_{1},\dots,v_{n}), \, \hat u=(u_{1}-2,\dots,u_{m}-2), \, \hat v=(v_{1}-2,\dots,v_{n}-2)$ and
\begin{equation}
\label{K}
K_{\lambda,\tilde\lambda}=4^{|\lambda|-|\tilde\lambda|}\frac{C_{\lambda}^0(d) \,C_{\lambda}^0(d-p-q-\frac12)}{C_{\tilde\lambda}^0(d)\, C_{\tilde\lambda}^0(d-p-q-\frac12)}\frac{I_{\tilde\lambda}(w(\lambda), z(\lambda),-1,h)}{C_{\lambda}^{-}(1)C_{\lambda}^+(2h-1)}.
\end{equation}
Here  $d=m-n$ is the superdimension, $h= d-\frac{1}{2}p-q,$
\begin{equation}
\label{C+}
C^{+}_{\lambda}(x)=\prod_{(ij)\in\lambda}\left(\lambda_{i}+j-(\lambda^{\prime}_{j}+i)+x\right),
\end{equation}
\begin{equation}
\label{C-}
C^{-}_{\lambda}(x)=\prod_{(ij)\in\lambda}\left(\lambda_{i}-j+(\lambda^{\prime}_{j}-i)+x\right),
\end{equation}
\begin{equation}
\label{C0}
C^{0}_{\lambda}(x)=\prod_{(ij)\in\lambda}\left(j-1-(i-1)+x\right),
\end{equation}
$I_{\lambda}(w, z,-1,h)$ is the specialisation of the deformed interpolation $BC$ polynomial $I_{\lambda}(w,z,k,h)$ (see Proposition 6.3 in \cite{SV2}), $w(\lambda)$ and $z(\lambda)$ are defined by (\ref{embedding}) above. We should note that we are using here more convenient variables $$u_i=x_i+x_i^{-1},\,v_j=y_j+y_j^{-1},$$ rather than
$u_i=\frac{1}{2}(x_i+x_i^{-1}-2), v_j=\frac{1}{2}(y_j+y_j^{-1}-2)$ used in \cite{SV2}.

\begin{thm}\label{factori}
Let $\lambda\in H(m,n)$ contains the $m\times n$ rectangle. Then the super Jacobi polynomials $SJ_{\lambda}(u,v,-1,p,q)$ can be expressed in terms of the usual Jacobi polynomials as 
\begin{equation}
\label{factor}
SJ_{\lambda}(u,v,-1,p,q)=
(-1)^{|\nu|}\prod_{i=1}^m\prod_{j=1}^n(u_{i}-v_{j})J_{\mu}(u,-1,p,q)J_{\nu}(v,-1, -p,-1-q),
\end{equation}
where 
\begin{equation}
\label{munu}\mu_{i}=\lambda_{i}-n,\: i=1,\dots,m,\quad \nu_{j}=\lambda^{\prime}_{j}-m,\; j=1,\dots,n.
\end{equation}
\end{thm}

\begin{proof}
To prove this we need the following factorisation formula for the deformed interpolation $BC$ polynomials.
We should mention that in the special case  $k=-1$ these polynomials are the particular case  of the factorial super Schur functions considered by Molev in \cite{Mo}, but in our particular case one can give a simple direct proof.

\begin{lemma}\label{factorisation1} If $\lambda \in H(m,n)$ contains the $m\times n$ rectangle, then we have the following formula for deformed interpolation $BC$ polynomials
$$
I_{\lambda}(w,z,-1, h)=
(-1)^{|\nu|}\prod_{i=1}^m\prod_{j=1}^n\left[(w_{i}+n+h-i)^2-(z_{j}-h-j+1)^2\right]\times
$$
$$
 I_{\mu}(w_{1},\dots,w_{m},-1,h+n)I_{\nu}(z_{1}-m,\dots,z_{n}-m,-1,1-h+m),
 $$
where $I_{\mu}(w,k,h)$ is the usual  interpolation $BC$ polynomial and $\mu,\,\nu$ are the same as in the theorem.
\end{lemma}

\begin{proof} Let us denote by $\hat I_{\lambda}$ the right hand side of the previous equality. 
According to proposition 6.3 from \cite{SV2} it is enough to prove that  $\hat I_{\lambda}$
satisfy the following properties:

1) $\hat I_{\lambda}$ is a polynomial in variables $(w_{i}+h+n-i)^2,\: i=1,\dots,m$ and  $(z_{j}-h+1-j)^2,\: j=1,\dots,n$;

2) the degree of this polynomial is $2|\lambda|$;

3) $\hat I_{\lambda}(w(\tilde\lambda),z(\tilde\lambda))=0$ if $\tilde\lambda\in H(m,n)$ and $\lambda\nsubseteq \tilde\lambda$;

4) $
\hat I_{\lambda}(w(\lambda),z(\lambda))=\prod_{(i,j)\in\lambda}(1+\lambda_{i}-j+\lambda^{\prime}_{j}-i)(2h-1+\lambda_{i}+j-\lambda^{\prime}_{j}-i).
$

First two statements are obvious from the explicit form of $\hat I_{\lambda}$. The fourth statement can be checked directly. Let us prove the third property. 

Let us suppose that $\lambda\nsubseteq \tilde\lambda$. Consider  two possible cases depending on whether $\tilde\lambda$ contains the $m\times n$ rectangle or not. 
In the first case we have $\mu\nsubseteq\tilde\mu$ or  $\nu\nsubseteq\tilde\nu$. Therefore by definition  of the interpolation $BC$ polynomials \cite{OO} $$I_{\mu}(\tilde\mu,-1,h+n)I_{\nu}(\tilde\nu,-1,1-h+m)=0.$$
In the second case consider  the box $(i,n),\: 1\le i\le m$ such that $(i,n)\notin\tilde\lambda$, but $(i-1,n)\in\tilde\lambda$ (if i=1, we require only first condition). Note that from our assumptions on $\lambda$ it follows  that such a box does exist. Then we have $(w_{i}+n+h-i)+(z_{n}-h-n+1)=\tilde\mu_{i}+n-i+\tilde\nu_{n}-n+1=0+n-i+(i-1)-n+1=0$, which means that $\hat I_{\lambda} (w(\tilde\lambda),z(\tilde\lambda))=0$. Lemma is proved.
\end{proof}

Now we can prove the theorem \ref{factori}.
First we note that in Okounkov's formula we can always assume that the diagram $\tilde\lambda$
contains the $m\times n$ rectangle. Indeed, otherwise $K_{\lambda,\tilde\lambda}=0$ since 
in that case one can easily see that $$\frac{C_{\lambda}^0(d)}{C_{\tilde\lambda}^0(d)}= 0.$$
Therefore by lemma \ref{factorisation1}
$$
I_{\tilde\lambda}(w(\lambda), z(\lambda),-1,h)=(-1)^{|\tilde\nu|}\prod_{i=1}^m\prod_{j=1}^n(\lambda_{i}+\lambda^{\prime}_{j}-i-j+1)(\lambda_{i}-\lambda^{\prime}_{j}+j-i+2h-1)\times
$$
$$
I_{\tilde\mu}(\mu,-1, h+n)I_{\tilde\nu}(\nu,-1,1- h+m),
$$
where $\tilde\mu, \tilde\nu$ are defined by $\tilde\lambda$ as in (\ref{munu}).
We rewrite now the coefficient $K_{\lambda,\tilde\lambda}$ in terms of the diagrams $\mu$ and $\nu.$
We have
$$
C_{\lambda}^{-}(1)=\prod_{i=1}^m\prod_{j=1}^n(\lambda_{i}+\lambda^{\prime}_{j}-i-j+1)C_{\mu}^{-}(1)C_{\nu^{\prime}}^{-}(1),
$$
$$
C_{\lambda}^{+}(2h-1)=\prod_{i=1}^m\prod_{j=1}^n(\lambda_{i}-\lambda^{\prime}_{j}+j-i+2h-1)C_{\mu}^{+}(2(h+n)-1)C_{\nu^{\prime}}^{+}(2(h-m)-1),
$$
$$
\frac{C_{\lambda}^0(m-n)}{C_{\tilde\lambda}^0(m-n)}=\frac{C_{\mu}^0(m)}{C_{\tilde\mu}^0(m)}\frac{C_{\nu^{\prime}}^0(-n)}{C_{\tilde\nu^{\prime}}^0(-n)},
$$
$$
\frac{C_{\lambda}^0(m-n-p-q-\frac12)}{C_{\tilde\lambda}^0(m-n-p-q-\frac12)}=\frac{C_{\mu}^0(m-p-q-\frac12)}{C_{\tilde\mu}^0(m-p-q-\frac12)}\frac{C_{\nu^{\prime}}^0(-n-p-q-\frac12)}{C_{\tilde\nu^{\prime}}^0(-n-p-q-\frac12)}.
$$
Now we use the Berele-Regev factorisation formula \cite{BR} for super Schur polynomials
$$
SP_{\tilde\lambda}(u,v,-1)=(-1)^{|\tilde\nu|}\prod_{i=1}^m\prod_{j=1}^n(u_{i}-v_{j})P_{\tilde\mu}(u)P_{\tilde\nu}(v),
$$
where $P_{\lambda}(u)$ are usual Schur polynomials.
Comparing (\ref{okoun}) and (\ref{factor}) and using  Okounkov's formula for usual Jacobi polynomials \cite{OO} we see that in order to prove the theorem we need to show that
$$
J_{\nu}(v_{1},\dots,v_{n},-1, -p,-1-q)=(-1)^{|\nu|}\sum_{\tilde\nu\subseteq\nu}D_{\nu,\tilde\nu}P_{\tilde\nu}(\hat v),
$$
where
$$
D_{\nu,\tilde\nu}=4^{|\nu|-|\tilde\nu|}\frac{C_{\nu^{\prime}}^0(-n)}{C_{\tilde\nu^{\prime}}^0(-n)}\frac{C_{\nu^{\prime}}^0(-n-p-q-\frac12)}{C_{\tilde\nu^{\prime}}^0(-n-p-q-\frac12)}\frac{I_{\tilde\nu}(\nu,-1,1- h+m)}{C_{\nu^{\prime}}^{-}(1)C_{\nu^{\prime}}^{+}(2(h-m)-1)}.
$$
Since
$$
C_{\nu}^{0}(x)=(-1)^{|\nu|}C_{\nu^{\prime}}^{0}(-x),\;C_{\nu}^{-}(x)=C_{\nu^{\prime}}^{-}(-x),\:C_{\nu}^{+}(x)=(-1)^{|\lambda|}C_{\lambda}^{+}(-x),
$$
we have
$$
(-1)^{|\nu|}D_{\nu,\tilde\nu}=4^{|\nu|-|\tilde\nu|}\frac{C_{\nu}^0(n)}{C_{\tilde\nu}^0(n)}\frac{C_{\nu}^0(n+p+q+\frac12)}{C_{\tilde\nu}^0(n+p+q+\frac12)}\frac{I_{\tilde\nu}(\nu,-1,1- h+m)}{C_{\nu}^{-}(1)C_{\nu}^{+}(1-2(h-m))}
$$
and the proof now follows from Okounkov's formula. 
\end{proof}

Let's come to the proof of the main theorem \ref{main}.
We need the following result. Denote the right hand side of the formula (\ref{WeylJ}) as $RHS.$

\begin{proposition} \label{propos} The right hand side of the formula (\ref{WeylJ}) can be written in terms of super Schur functions as
$$
RHS= \sum_{\mu \subseteq \lambda}C_{\mu}(d,p,q)SP_{\mu}(u_{1}-2,\dots,u_{m}-2,v_{1}-2,\dots,v_{n}-2),
$$
where the coefficients $C_{\mu}(d,p,q)$ are some rational functions of $\mu,d,p,q.$
\end{proposition}

To prove this let us denote $r=i(\lambda),\; s=j(\lambda)$ and
expand every polynomial $f_{l_{i}},\; 1\le i\le r$ in terms of the powers of $u_{i}-2$ and every polynomial $g_{k_{j}},\; 1\le j\le s$ in terms of the powers of $v_{j}-2$.  Then RHS is the sum of terms $$
\sum_{w\in S_{m}\times S_{n}}w\left[\Pi_{\lambda}(u,v)\frac{(u_{1}-2)^{\tilde l_{1}}\dots(u_{r}-2)^{\tilde l_{r}}(v_{1}-2)^{\tilde k_{1}}\dots(v_{r}-2)^{\tilde l_{s}}F_{\lambda}G_{\lambda}}{\Delta(u)\Delta(v)}\right]
$$
with some constant factors depending on $d,p,q,$ 
where 
$$F_{\lambda}=(u_{r+1}-2)^{m-r-1}\dots(u_{m}-2)^0,\: G_{\lambda}=(v_{s+1}-2)^{m-s-1}\dots(v_{n}-2)^0$$
and $ T=(\tilde l_{1},\dots,\tilde l_{r},\tilde k_{1},\dots\tilde k_{s})$ satisfy the  conditions  $ 0\le\tilde l_{i}\le l_{i}\:, 1\le i\le i(\lambda)$ and  $ 0\le\tilde k_{j}\le k_{j}\:, 1\le j\le j(\lambda)$.  Since $\pi_{\lambda}$ is symmetric with respect to $u_{1},\dots,u_{r}$ and $v_{1},\dots,v_{s}$ we may assume that both $\tilde l_{i}$ and $\tilde k_{j}$ are pairwise different. Now the proposition follows from theorem \ref{SS} because of the following
\begin{lemma} Let $a_{1}>\dots>a_{n}$ sequence of nonnegative integers  and $b_{1},\dots,b_{n}$ sequence of pairwise different nonnegative integers such that $a_{i}\ge b_{i}, i=1,\dots,n$.  Let us reorder  sequence $\{b_{i}\}$ in decreasing order $b_{1}^{\prime}>b_{2}^{\prime}>\dots>b_{n}^{\prime}$. Then for any $1\le i\le n$ we have $b_{i}^{\prime}\le a_{i}$.
\end{lemma}
\begin{proof}
We will prove lemma by induction on $n$. The case $n=1$ is obvious.
Assume that lemma is true for some $n$ and consider the sequences $a_{1}>\dots>a_{n}>a_{n+1},\; b_{1},\dots,b_{n},b_{n+1},$ which satisfy the conditions of the lemma and the corresponding sequence  $b_{1}^{\prime}>b_{2}^{\prime}>\dots>b_{n}^{\prime}>b_{n+1}^{\prime}$. Let $c=\min\{b_{1},\dots,b_{n+1}\}$. If $c= b_{n+1}.$ We can apply inductive assumption to $a_{1},\dots,a_{n},\:b_{1},\dots,b_{n}$. If $c=b_{i}\ne b_{n+1}$ then we can apply inductive assumption to $a_{1},\dots,a_{n}$ and $b_{1},\dots,b_{i-1},b_{n+1},b_{i+1},\dots,b_{n}$. Lemma is proved.
\end{proof}

Let us finish the proof of the theorem \ref{main}. The proof is by induction on $m+n.$ 
First note that if $\lambda$ contains the $m\times n$ rectangle then the theorem
follows from Theorem \ref{factori} and Proposition 7.1 from Okounkov-Olshanski \cite{OO}.
In particular, this is always true when either $m$ or $n$ is zero.

Assume now that both $m$ and $n$ are positive and $\lambda$ does not contain the $m\times n$ rectangle. Make the substitution $u_m=v_n=t$ in both sides of the formula (\ref{WeylJ}).
The result is independent on $t$ and reduces to the case with smaller number of variables $u_1, \dots, u_{m-1}$ and $v_1, \dots v_{n-1}$: for the left hand side this follows from Okounkov's formula, for the right hand side from proposition \ref{propos}. Thus by induction we have that the difference between
the right hand side and  left hand side of the formula (\ref{WeylJ}) is divisible by the product
$\prod_{i=1}^m\prod_{j=1}^n (u_i-v_j).$ Note that both sides are linear combinations of the super Schur polynomials $SP_{\tilde\lambda}$ with $\tilde\lambda \subseteq \lambda$, which follows from Okounkov's formula and proposition \ref{propos}. However the ideal generated by $\prod_{i=1}^m\prod_{j=1}^n (u_i-v_j)$ is known to be linearly spanned by the super Schur polynomials $\mathcal SP_{\mu}$ with $\mu$ containing the $m\times n$ rectangle. This means that the difference is actually zero since by assumption $\lambda$ does not contain it. This completes the proof of Theorem \ref{main}.

As a corollary we have one of our main results. 
\begin{thm} \label{ejtheor}
The limit of the super Jacobi polynomials $SJ_{\lambda}(u,v,-1,p,q)$ as $(p,q) \rightarrow (-1,0)$ is well defined and coincides up to a constant factor with the Euler supercharacter for the Lie superalgebra $\mathfrak{osp}(2m+1,2n).$:
\begin{equation}
\label{ej}
\mathcal SJ_{\lambda}(u,v; -1,-1,0) = 2^{i(\lambda)-m}E_{\lambda}(u,v).
\end{equation}
\end{thm}
The proof follows from comparison of formulas (\ref{Euler}) and (\ref{WeylJ}).

\section{Pieri formula}

One of the problems with the special value $k=-1$ is that the spectrum of the corresponding ring of quantum integrals is not simple. This means that we need more information to characterize the specialized super Jacobi polynomials in this case. In this section we derive the Pieri formula for super Jacobi  polynomials in the case when $k=-1$ and show that it allows to characterize them uniquely.

This formula can be deduced from the Pieri formula \cite{SV2} for general $k$ by taking the limit $k\rightarrow -1$. However the corresponding calculations are quite long, so we will do this in a different  way.

Let us introduce some notations.
Let $\mathcal P_{m,n}$ be the set of all partitions $\lambda$ from the fat $(m,n)$ hook $H(m,n)$. Let us call a box $\Box=(i,j)$ of Young diagram $\lambda$
{\it special} if $$i-j=d,$$
where $d=m-n$ as before is superdimension. We will write $\mu\sim\lambda $ if Young diagram
$\mu$ can obtained from $\lambda$ by removing or adding one box.

Introduce the following functions:
\begin{equation}
\label{ad}
a_d(\mu,\lambda)=\begin{cases}0,\:\: \mu=\lambda\setminus\Box\,\, \text {for special} \,\,  \Box \cr 1,\: \:\,\, \text{otherwise},\end{cases}
\end{equation}
\begin{equation}
\label{bd}
b_d(\lambda)=\begin{cases}-1,\:\: \text{if $\lambda$ has a removable special box}\cr 1,\: \:\,\, \text{if there is a special box which can be added to $\lambda$}\cr 0,\:\:\,\, \text{otherwise.}\end{cases}
\end{equation}
\begin{thm}\label{Pieri} Let $\mathcal SJ_{\lambda}(u,v; -1,-1,0)=SJ_{\lambda}(u,v)$ be  specialised super Jacobi 
polynomials,  then the following Pieri formula holds
$$
\left(\sum_{i=1}^m u_{i}-\sum_{j=1}^n v_{j}+1\right)\mathcal SJ_{\lambda}(u,v)=$$
\begin{equation}
\label{pieri}
\sum_{\mu \sim \lambda, \, \mu\in \mathcal P_{m,n}}a_d(\mu,\lambda)\mathcal SJ_{\mu}(u,v)
+b_d(\lambda)\mathcal SJ_{\lambda}(u,v)
\end{equation}

\end{thm} 
\begin{proof}  We need the following Pieri formula for the Jacobi symmetric functions $\mathcal J_{\lambda}(u,k,p,q,h) \in \Lambda,$ where $\Lambda$ is the algebra of symmetric functions \cite{Mac}, $h$ is an additional parameter (see \cite{SV2}). First note that the limit of $\mathcal J_{\lambda}(u,k,p,q,h)$ when $k \rightarrow -1$ for generic $p,q$ is well defined as it follows from Okounkov's formula. We denote this limit $\mathcal J_{\lambda}.$ By $\lambda\pm\varepsilon_{i}$ we denote  the sets $(\lambda_1, \dots,\lambda_{i-1}, \lambda_i\pm1,\lambda_{i+1}, \dots)$ respectively. 

\begin{lemma}\label{Pieriinf} For $k=-1$ and generic $p,q$ the Jacobi symmetric functions satisfy the following Pieri formula
\begin{equation}
\label{pieriinf}
p_{1}\mathcal J_{\lambda}=
\sum_{i:\lambda+\varepsilon_{i}\in \mathcal P_{m,n}}\mathcal J_{\lambda+\varepsilon_{i}}+
\left(\sum_{i=1}^{l(\lambda)}a(\lambda_{i}+d-i)\right)\mathcal J_{\lambda}
\end{equation}
$$
+\left(p+\frac{p(p+2q+1)}{2h-2l(\lambda)+1}\right)\mathcal J_{\lambda}+\sum_{i:\lambda-\varepsilon_{i}\in \mathcal P_{m,n}}b(\lambda_{i}+d-i)\mathcal J_{\lambda-\varepsilon_{i}}
$$
where $p_1=u_1+u_2+\dots,$ $d=h+\frac{1}{2}p+q$ and  
\begin{equation}
\label{3a}
a(l)=-\frac{2p(p+2q+1)}{(2l-p-2q-1)(2l-p-2q+1)},
\end{equation}
\begin{equation}
\label{3b}
 b(l)=\frac{2l(2l-2q-1)(2l-2p-2q-1)(2l-2p-4q-2)}{(2l-p-2q)(2l-p-2q-1)^2(2l-p-2q-2)}.
\end{equation}
\end{lemma}
To prove the lemma we use the following formula for the Jacobi polynomials with $k=-1$ in $m$ variables (see Proposition 7.1 in Okounkov-Olshanski \cite{OO}):
$$
J_{\lambda}(u,-1,p,q)=\frac{1}{\Delta(u)}\sum_{w\in S_{m}}
\varepsilon(w)w\left[f_{\lambda_{1}+m-1}(u_{1})f_{\lambda_{2}+m-2}(u_{2})\dots f_{\lambda_{m}}(u_{m})\right],
$$ 
where as before $f_{l}(z)=f_l(z,p,q)$ are the classical normalized Jacobi polynomials in one variable with parameters $p,q$. They satisfy the following  three-term recurrence relation (see e.g. \cite{Erd}):
$$
zf_{l}(z)=f_{l+1}(z)+a(l)f_{l}(z)+b(l)f_{l-1}(z)
$$
with $a(l)$, $b(l)$ given by (\ref{3a}),(\ref{3b}).
Therefore we have
$$
\left(\sum_{i=1}^mu_{i}\right)J_{\lambda}(u)=\frac{1}{\Delta(u)}\sum_{i=1}^m\sum_{w\in S_{m}}\varepsilon(w)w\left[u_{i}f_{\lambda_{1}+m-1}(u_{1})\dots f_{\lambda_{m}}(u_{m})\right]
$$
$$
=\sum_{i:\lambda+\varepsilon_{i}\in \mathcal P_{m,n}}J_{\lambda+\varepsilon_{i}}(u)+\left(\sum_{i=1}^{l(\lambda)}a(\lambda_{i}+m-i)+\sum_{i=l(\lambda)+1}^{m}a(m-i)\right)J_{\lambda}(u)
$$
$$+\sum_{i:\lambda-\varepsilon_{i}\in \mathcal P_{m,n}}b(\lambda_{i}+m-i)J_{\lambda-\varepsilon_{i}}(u),
$$
where $J_{\lambda}(u) = J_{\lambda}(u,-1,p,q)$ and $l(\lambda)$ is the number of non-zero parts in partition $\lambda.$ Since
$$
a(x)=\frac{p(p+2q+1)}{2x-p-2q+1}-\frac{p(p+2q+1)}{2x-p-2q-1}
$$
we have
$$
\sum_{i=l(\lambda)+1}^{m}a(m-i)=p+\frac{p(p+2q+1)}{2m-2l(\lambda)-p-2q-1}.
$$
Comparing this with (\ref{pieriinf}) we see that this formula is true after the natural homomorphism $\Lambda \rightarrow \Lambda_m$, $\Lambda_m$ is the algebra of symmetric polynomials of $m$ variables, if we specialize $d$ to $m.$ Since this is valid for all $m,$ the lemma follows.

The super Jacobi polynomials $SJ_{\lambda}(u,v,-1,p,q)$ are defined as the image of Jacobi symmetric functions $\mathcal J_{\lambda}$ under  the homomorphism  
$$
\varphi : \Lambda\rightarrow\Lambda_{m,n},\quad\varphi(p_{l})=\sum_{i=1}^mu^l_{i}-\sum_{j=1}^nv_{j}^l,
$$ 
where $d$ is specialized to $m-n$ and $\Lambda_{m,n}$ is the algebra of supersymmetric polynomials (see e.g. \cite{Mac}).
Computing the limits when $(p,q)\rightarrow(-1,0)$
$$
\lim_{(p,q)\rightarrow(-1,0)}a(l)=\delta(l+1)-\delta(l),\;\; \lim_{(p,q)\rightarrow(-1,0)}b(l)=1-\delta(l)
$$
$$
\lim_{(p,q)\rightarrow(-1,0)}\left(p+\frac{p(p+2q+1)}{2d-2l(\lambda)-p-2q-1}\right)=-1+\delta(d-l(\lambda)),
$$
we have the following formula
$$
\left(\sum_{i=1}^m u_{i}-\sum_{j=1}^n v_{j}\right)\mathcal SJ_{\lambda}(u,v)=$$
\begin{equation}
\label{pieriold}
\sum_{i:\lambda+\varepsilon_{i}\in \mathcal P_{m,n}}\mathcal SJ_{\lambda+\varepsilon_{i}}(u,v)
+\sum_{i:\lambda-\varepsilon_{i}\in \mathcal P_{m,n}}\left[1-\delta(\lambda_{i}-i+d)\right]\mathcal SJ_{\lambda-\varepsilon_{i}}(u,v)
\end{equation}
$$
+\sum_{i=1}^{\l(\lambda)}\left[\delta(\lambda_{i}-i+d+1)-\delta(\lambda_{i}-i+d)\right]\mathcal SJ_{\lambda}(u,v)+
\left[\delta(d-l(\lambda))-1\right]SJ_{\lambda}(u,v),
$$
where $d=m-n\:$ and 
$$
\delta(x)=\begin{cases}1,\:\: x=0\cr 0,\: \: x\ne0.\end{cases}.
$$ 
One can check that it is equivalent to the
Pieri formula (\ref{pieri}).
\end{proof}

\begin{remark} The form of the factor on the left hand side of Pieri formula (\ref{pieri}), which was chosen by convenience, has a clear representation-theoretic meaning: $\sum_{i=1}^m u_{i}-\sum_{j=1}^n v_{j}+1$ is the supercharacter of the standard representation of $\mathfrak{osp}(2m+1,2n).$
\end{remark}

\begin{remark} One can possibly use this for the alternative proof of the main theorem. For this one should only prove that the Euler supercharacters satisfy the corresponding version of the Pieri formula. This is related to the translation functors in representation theory (see e.g. \cite{Brun}).
\end{remark}

\section{Jacobi--Trudy formula for Euler supercharacters}

In this section we give one more formula for specialized super Jacobi polynomials (and hence for the Euler supercharacters) of Jacobi--Trudy type.

We start with the Jacobi--Trudy formula for Jacobi symmetric functions following \cite{SV4}. 
Let $\mathcal J_{\lambda} \in \Lambda$ be the Jacobi symmetric functions with $k=-1$ and generic $p$ and $q$ (see \cite{SV2}). They depend also on the additional parameter $d$ replacing the dimension of the space. Let $h_i=\mathcal J_{\lambda}$ for the partition $\lambda=(i)$ consisting of one part for positive $i$ and $h_i\equiv 0$ if $i<0.$
Define recursively the sequence $h^{(r)}_i\in \Lambda, i\in\mathbb Z$ by the relation
\begin{equation}
\label{relahinf}
h_i^{(r+1)}=h_{i+1}^{(r)}+ a(i+d-1)h_i^{(r)}+b(i+d-1)h_{i-1}^{(r)}
\end{equation}
for $r=0,1,\dots$ with initial data $h_i^{(0)}=h_i$ and 
$$
a(x)=-\frac{2p(p+2q+1)}{(2x-p-2q-1)(2x-p-2q+1)},
$$
$$
 b(x)=\frac{2x(2x-2q-1)(2x-2p-2q-1)(2x-2p-4q-2)}{(2x-p-2q)(2x-p-2q-1)^2(2x-p-2q-2)}.
$$

\begin{thm} \cite{SV4} The Jacobi symmetric functions with $k=-1$ have the following Jacobi-Trudy representation:
\begin{equation}
\label{JTinf}
\mathcal J_{\lambda}=\left|\begin{array}{cccc}
   h_{\lambda_{1}}&h^{(1)}_{\lambda_{1}}& \ldots &h^{(l-1)}_{\lambda_{1}}\\
  h_{\lambda_{2}-1}&h^{(1)}_{\lambda_{2}-1}& \ldots &h^{(l-1)}_{\lambda_{2}-1}\\
\vdots&\vdots&\ddots&\vdots\\
  h_{\lambda_{l}-l+1}&h^{(1)}_{\lambda_{l}-l+1}& \ldots &h^{(l-1)}_{\lambda_{l}-l+1}\\
 \end{array}\right|,
 \end{equation}
where $l=l(\lambda)$. 
\end{thm}

Taking the limit $p\rightarrow -1,\, q\rightarrow 0$ and the homomorphism $\phi_{m,n}:\Lambda\rightarrow \Lambda_{m,n}$ we have the following 

\begin{corollary} The specialized super Jacobi polynomials $\mathcal SJ_{\lambda}$
satisfy the following Jacobi--Trudy formula
\begin{equation}
\label{JTinfsup}
\mathcal SJ_{\lambda}(u,v;-1,-1,0)=\left|\begin{array}{cccc}
   h_{\lambda_{1}}&h^{(1)}_{\lambda_{1}}& \ldots &h^{(l-1)}_{\lambda_{1}}\\
  h_{\lambda_{2}-1}&h^{(1)}_{\lambda_{2}-1}& \ldots &h^{(l-1)}_{\lambda_{2}-1}\\
\vdots&\vdots&\ddots&\vdots\\
  h_{\lambda_{l}-l+1}&h^{(1)}_{\lambda_{l}-l+1}& \ldots &h^{(l-1)}_{\lambda_{l}-l+1}\\
 \end{array}\right|,
 \end{equation}
where $\lambda \in H_{m,n},\,d=m-n$  and $h^{(r)}_i$ are defined recursively by (\ref{relahinf}) with
$$
a(x)=\delta(x+1)-\delta(x),\quad b(x)=1-\delta(x)
$$
and $h_i^{(0)}=h_i=\mathcal SJ_{\lambda}(u,v;-1,-1,0)$ for $\lambda=(i)$ for positive $i$ and $h_i^{(0)}=h_i\equiv 0$ if $i<0.$
\end{corollary}

\section{Euler supercharacters for different choice of Borel subalgebra}

It turns out that the relation with super Jacobi polynomials can be made more direct if we choose, following to Gruson and Serganova \cite{GS}, a different Borel subalgebra and suitable parabolic subalgebras. 

We consider again the case $\frak{g}=\frak{osp}(2m+1,2n)$, assuming for convenience at the beginning that $m\ge n$. Choose the following set of simple roots
$$
B=\{\varepsilon_{1}-\varepsilon_{2},\dots,\varepsilon_{m-n+1}-\delta_{1},\delta_{1}-\varepsilon_{m-n+2},\varepsilon_{m-n+2}-\delta_2,\dots,\varepsilon_{m}-\delta_{n}, \delta_{n}\}
$$
This choice is special since this set contains the maximal possible number of isotropic roots.
The corresponding set of even positive roots is
$$
R^+_{0}=\{\varepsilon_{i}\pm\varepsilon_{j},\,\,i <j,\,\,\,\,\delta_{p}\pm\delta_{q},\,p<q,\,2\delta_{p}\}
$$
with the half-sum
$$
\rho_{0}=\frac12\sum_{\alpha\in R^+_0}\alpha=(m-\frac12)\varepsilon_{1}+(m-\frac32)\varepsilon_{2}\dots+\frac12\varepsilon_{m}+n\delta_{1}+(n-1)\delta_{2}+\dots+\delta_{n}.
$$
The set of positive odd roots is
$$
R^+_{1}=\{\varepsilon_{i}\pm\delta_{j},\,\,i -j\le m-n,\,\,\,\,\delta_{j}\pm\varepsilon_{i},\,i-j>m-n,\,\delta_{j},\, j=1,\dots, n\}
$$
with the half-sum
$$
\rho_{1}=n(\varepsilon_{1}+\dots+\varepsilon_{m-n+1})+(n-1)\varepsilon_{m-n+2}+\dots +\varepsilon_{m}+(n-\frac12)\delta_{1}+\dots+\frac12\delta_{n},
$$
so
$$
\rho=\rho_{0}-\rho_{1}=\sum_{i=1}^{m-n}(m-n-i+\frac12)\varepsilon_{i}-\frac12\sum_{i=m-n+1}^m\varepsilon_{i}+\frac12\sum_{j=1}^n\delta_{j}.
$$

The following lemma gives a description of the highest weights with respect to our choice of simple roots $B$ in terms of partitions. We need this to establish the relation with super Jacobi polynomials.

\begin{lemma}\label{Young} The weight 
$$
\chi=a_{1}\varepsilon_{1}+a_{2}\varepsilon_{2}+\dots+a_m\varepsilon_{m}+b_{1}\delta_{1}+b_{2}\delta_{2}+\dots+b_n\delta_{n}
$$
is a highest weight of an irreducible finite dimensional $\frak{osp}(2m+1,2n)$-module   if and only if there exists a partition $\lambda \in H(m,n)$  from the fat hook such that
\begin{equation}
\label{weight}
\chi+\rho= \sum_{i=1}^{i(\lambda)}(\lambda_{i}+d-i+\frac12)\varepsilon_i - \frac12\sum_{i>i(\lambda)}\varepsilon_i+ 
\sum_{j=1}^{j(\lambda)}(\lambda^{\prime}_{j}-d-j+\frac12) \delta_j + \frac12\sum_{j>j(\lambda)}\delta_j.
\end{equation}
\end{lemma}
The proof is a geometric reformulation of the conditions on the highest weights given in \cite{GS}, Corollary 3.
 
Let $\lambda$ be a partition from $H(m,n)$ and $\chi$ be the corresponding highest weight. Consider the maximal parabolic subalgebra $\mathfrak p=\mathfrak p(\lambda)$  such that $\chi$ can be extended to $\mathfrak p$ as a one dimensional representation. One can check that in that case  the corresponding roots are
$$
R_{\mathfrak p}=\{\pm\varepsilon_{i}\pm\varepsilon_{j},\pm\varepsilon_{i},\,\,\pm\delta_{p}\pm\delta_{q},\pm\delta_{p},\,\,\pm2\delta_{p},\,  \pm \varepsilon_i \pm \delta_p\}
$$ 
with
$\,i(\lambda)\le i,j\le m,\,\,i\ne j,\,\,\,\,j(\lambda)\le p,q\le n,\,\,p\ne q\,.$ They correspond to the subalgebra $\mathfrak{osp}(2l+1,2l),\, l = m-i(\lambda)=n-j(\lambda).$

Let $E^{\mathfrak p}(\chi)$ be the corresponding Euler supercharacter given by (\ref{Vera}).
Let $\pi_{\lambda}$ be the same as in section 2 and define $t(\lambda)$ as the number of pairs $(i,j) \in \pi(\lambda)$ with $i-j>d=m-n.$ Similarly let $s(\lambda)$ be the number of pairs $(i,j) \in \lambda$ with $i-j>d.$

\begin{thm}
Let $\mathfrak p$ and $\chi$ be the parabolic subalgebra and its one-dimensional representation determined by a partition $\lambda \in H(m,n)$, then the specialized Jacobi polynomial $SP_{\lambda}(u,v,-1,-1,0)$ coincides up to a sign with the corresponding Euler supercharacter: 
\begin{equation}
\label{main2}
SP_{\lambda}(u,v,-1,-1,0) = (-1)^{s(\lambda)} E^{\mathfrak p}(\chi).
\end{equation}
\end{thm}

\begin{proof}
According to general formula (\ref{Vera})
$$E^{\mathfrak p}(\chi) = \sum_{w \in W_0} \, w\left(\frac{\prod_{\alpha \in R_1^+}(e^{\alpha/2} - e^{-\alpha/2}) \, e^{\chi+\rho}}{\prod_{\alpha \in R_0^+}(e^{\alpha/2} - e^{-\alpha/2})\prod_{\alpha \in R_{\mathfrak p} \cap R_1^+} (1-e^{-\alpha})}\right)=
$$
$$
\sum_{w \in W_0} \, w\left(\frac{ \prod_{\alpha \notin R_{\mathfrak p} \cap R_1^+} (e^{\alpha/2} - e^{-\alpha/2})\, e^{\chi+\rho+\tau}}{\prod_{\alpha \in R_0^+}(e^{\alpha/2} - e^{-\alpha/2})}\right),
$$
where $$\tau = \frac{1}{2} \sum_{\alpha \in R_{\mathfrak p} \cap R_1^+} \alpha.$$
A simple calculation shows that
$$
 \prod_{\alpha \notin R_{\mathfrak p} \cap R_1^+} (e^{\alpha/2} - e^{-\alpha/2})= (-1)^{t(\lambda)}\prod_{(i,j)\in \pi_{\lambda}} (u_i-v_j)\prod_{j\leq j(\lambda)} (y_j^{\frac{1}{2}}-y_j^{-\frac{1}{2}}),
$$
where $u_{i}=x_{i}+x_{i}^{-1}, \, x_{i}=e^{\varepsilon_{i}},\, i=1,\dots,m$, $ v_{j}=y_{j}+y_{j}^{-1}, \,  y_{j}=e^{\delta_{j}},\: j=1,\dots, n.$
One can check also that
$$
\chi+\rho+\tau= \sum_{i=1}^m (l_i+\frac{1}{2}) \varepsilon_i + \sum_{j=1}^{j(\lambda)} (k_j+\frac{1}{2}) \delta_j +\sum_{j>j(\lambda)}^n (n-j+1)  \delta_j,
$$
where $l_i,\, k_j$ are defined by (\ref{frob11}), (\ref{frob21}). Therefore we have
$$E^{\mathfrak p}(\chi)=(-1)^{t(\lambda)}\sum_{w \in W_0} \, w\left(\frac{\prod_{(i,j)\in \pi_{\lambda}} (u_i-v_j)\prod_{j\leq j(\lambda)} (y_j^{\frac{1}{2}}-y_j^{-\frac{1}{2}}) e^{\chi+\rho+\tau}}{\Delta(u)\Delta(v) \prod_{i=1}^m (x_i^{\frac12}-x_i^{-\frac12})\prod_{j=1}^n (y_j-y_j^{-1})}\right)=
$$
$$
(-1)^{t(\lambda)}\sum_{w \in W_0} \, w\left(\frac{{ \displaystyle\prod_{(i,j)\in \pi_{\lambda}} (u_i-v_j)} \, e^{\chi+\rho+\tau}}{\Delta(u)\Delta(v) { \displaystyle \prod_{i=1}^m (x_i^{\frac12}-x_i^{-\frac12})\prod_{j\leq j(\lambda)} (y_j^{\frac{1}{2}}+y_j^{-\frac{1}{2}})\prod_{j>j(\lambda)}^n (y_j-y_j^{-1})}}\right).
$$
Now we use the identities
$$
\sum_{w \in W_l} \, w\left(\frac{y_1^l y_2^{l-1}\dots y_l}{\Delta(v) \prod_{j=1}^l (y_j-y_j^{-1})}\right)=1=
\sum_{w \in W_l} \, w\left(\frac{y_1^{l-\frac12} y_2^{l-\frac 32}\dots y_l^{\frac 12}}{\Delta(v) \prod_{j=1}^l (y_j^{\frac 12}+y_j^{-\frac 12})}\right),
$$
where $W_l$ is the Weyl group of type $C_l \approx BC_l$, which is a semi-direct product of permutation group $S_l$ and $\mathbb Z_2^l.$ These identities follow from the Weyl (super)character formula for trivial representations of $\mathfrak{sp}(2l)$ and $\mathfrak{osp}(1,2l).$
This leads to 
$$E^{\mathfrak p}(\chi)=(-1)^{t(\lambda)}\sum_{w \in W_0} \, w\left(\frac{{ \displaystyle\prod_{(i,j)\in \pi_{\lambda}} (u_i-v_j)} \, x_1^{l_1+\frac12}\dots x_m^{l_m+\frac12} y_1^{k_1+\frac12}\dots y_n^{k_n+\frac12}}{\Delta(u)\Delta(v) { \displaystyle \prod_{i=1}^m (x_i^{\frac12}-x_i^{-\frac12})\prod_{j=1}^n (y_j^{\frac{1}{2}}+y_j^{-\frac{1}{2}})}}\right).
$$
Comparing this with Theorem \ref{main} and using the obvious relation
$$s(\lambda)=t(\lambda)+ b(\lambda), \,\, b(\lambda)= \sum_{i>m} \lambda_i$$
we have the claim.
\end{proof}

\section{The case of $\mathfrak{osp}(2m,2n)$}

In this section we present the results in the even orthosymplectic case $\mathfrak g=\mathfrak{osp}(2m,2n).$

It would be instructive to start with the special case $n=0$, i.e. with the usual orthogonal Lie algebra $\mathfrak{o}(2m).$ It has the root system of type $D_m$ with simple roots
$$
\{\varepsilon_{1}-\varepsilon_{2},\dots,\varepsilon_{m-1}-\varepsilon_{m}, \varepsilon_{m-1}+ \varepsilon_{m}\}.
$$
We have a symmetry $\varepsilon_{m} \rightarrow - \varepsilon_{m}$, corresponding to the (outer) automorphism $\theta$ of this Lie algebra. This automorphism appears in the description of the representations of the corresponding Lie group $O(2m),$ which consists of two connected components.

Recall (see \cite{FH}) that the highest weights 
$\mu= \mu_1\varepsilon_{1} + \dots + \mu_m \varepsilon_{m} $
of irreducible finite-dimensional representations of Lie algebra $\mathfrak{o}(2m)$ have the following form: all $\mu_i$ are either integer or half-integer and satisfy the inequalities:
$$\mu_1 \geq \mu_2 \geq \dots \geq \mu_{m-1} \geq |\mu_m|.$$
The half-integer $\mu$ correspond to the spinor representations and can not be extended to the representations of the orthogonal group $SO(2m),$ so we restrict ourselves by integer $\mu.$
Corresponding representation $V^{\mu}$ of $\mathfrak{o}(2m)$ can be extended to the full orthogonal group $O(2m)$ if and only if it is invariant under the automorphism $\theta,$ which is equivalent to $\mu_m=0.$ If $\mu_m \neq 0$ then one should consider the direct sum
$$W^{\mu}=V^{\mu}\oplus V^{\theta(\mu)},$$
which gives an irreducible representation of $O(2m).$

It is interesting that the sum of the corresponding Euler supercharacters appears also in the general orthosymplectic case $\mathfrak g=\mathfrak{osp}(2m,2n)$
as the limit of super Jacobi polynomials (see below), so these limits are natural to link with supergroup $OSP(2m,2n)$ rather than Lie superalgebra $\mathfrak{osp}(2m,2n).$

Now let us give precise formulation of the results.
We have $\mathfrak{g}= \mathfrak{g}_0 \oplus \mathfrak{g}_1$, where $\mathfrak{g}_{0}=\mathfrak{so}(2m)\oplus \mathfrak{sp}(2n)$ and
$\mathfrak{g}_{1}=V_{1}\otimes V_{2}$ where $V_{1}$ and $V_{2}$ are the
identical representations of $\mathfrak{so}(2m)$ and $\mathfrak{sp}(2n)$ respectively. Let
$\pm\varepsilon_{1},\dots,\pm\varepsilon_{m},\, \pm\delta_1, \dots, \pm\delta_n$ be the non-zero weights of the
identical representation of  $\mathfrak{g}.$
The root system of $\mathfrak{osp}(2m,2n)$ consists of
$$R_{0}=\{\pm\varepsilon_{i}\pm\varepsilon_{j},\: \, \pm\delta_{p}\pm\delta_{q},\:\pm2\delta_{p}, \,
\: i\ne j,\,\: 1\le i,j\le m\ , p\ne q,\, 1\le p,q\le n\},
$$
$$
R_{1}=R_{iso}=\{\pm\varepsilon_{i}\pm\delta_{p}\},
$$
where $R_{0}, R_{1}$ and $R_{iso}$ are even, odd and isotropic parts respectively.
The Weyl group $W_0=\left(S_{m}\ltimes\mathbb{Z}_{2}^{(m-1)}\right)\times\left( S_{n}\ltimes\mathbb{Z}_{2}^n\right)$ acts on the weights by separately  permuting $\varepsilon_{i},\; j=1,\dots,m$ and
$\delta_{p},\; p=1,\dots,n$ and changing their signs such that the total number of signs of $\varepsilon_{i}$ is even.
A distinguished system of simple roots can be chosen as
$$
B=\{\delta_{1}-\delta_{2},\dots,\delta_{n-1}-\delta_{n},\delta_{n}-\varepsilon_{1},\varepsilon_{1}-\varepsilon_{2},\dots,\varepsilon_{m-1}-\varepsilon_{m}, \varepsilon_{m-1}+ \varepsilon_{m}\}.
$$
We have again a symmetry $\varepsilon_{m} \rightarrow - \varepsilon_{m}$, corresponding to the automorphism  of Lie superalgebra $\mathfrak{osp}(2m,2n)$, which also denote $\theta.$ It acts also in a natural way on the Grothendieck ring and supercharacters.

Consider the parabolic subalgebra $\mathfrak p$ with 
$$R_{\mathfrak p} = \{\varepsilon_{i}-\varepsilon_{j}, \, \delta_{p}-\delta_{q},\,\, \pm(\varepsilon_{i}- \delta_{p})\},$$ which is isomorphic to the sum of the Lie superalgebra $\mathfrak{gl}(m,n)$ and some nilpotent Lie superalgebra. 
For any finite-dimensional representation $M$ of $\mathfrak p$  the corresponding Euler supercharacter $E^{\mathfrak p}(M)$ is given by (\ref{Vera}), which in this case 
has the form 
\begin{equation}
\label{Vera2}
E^{\mathfrak p}(M) = \sum_{w \in W_0} w \left( \frac{\displaystyle\prod_{(i,j)\in \Pi_{m,n}}(1-x_i^{-1}y_j^{-1}) x_1^{m-1}\dots x_m^{0} y_1^{n}\dots y_n^{1} sch M}{\Delta(u)\Delta(v)\prod_{j=1}^n (y_j^{}-y_j^{-1})} \right).
\end{equation}
Here $\Pi_{m,n}$ is the rectangle of the size $m\times n$, $\Delta(u) = \prod_{i<j}^m(u_i-u_j),\, \Delta(v)=\prod_{i<j}^n(v_i-v_j)$ and $u_{i}=x_{i}+x_{i}^{-1}, \, x_{i}=e^{\varepsilon_{i}}, \, v_{j}=y_{j}+y_{j}^{-1},\, y_{j}=e^{\delta_{j}}$ as before.

\begin{proposition} \label{Triv1} For the trivial even representation $M$ we have
\begin{equation}
\label{triv1}
E^{\mathfrak p}(M) = 2^{\min(m-1,n)}.
\end{equation}
\end{proposition}

\begin{proof} The proof is similar to Proposition 3.1 and is based on the formula
$$\prod_{i,j}(1-x_i^{-1}y_j^{-1})= \sum (-1)^{|\lambda|} S_{\lambda}(x_1^{-1},\dots, x_m^{-1})S_{\lambda'}(y_1^{-1},\dots, y_n^{-1}),$$
where $S_{\lambda}$ is the Schur polynomial  and the sum is over all the Young diagrams $\lambda$, which are contained in the $(m\times n)$ rectangle. 
Replacing  in the formula (\ref{Vera2}) as before the product
$\prod_{i,j}(1-x_i^{-1}y_j^{-1})$ by the sum
 $$\sum_{\lambda \subseteq \Pi_{m,n}} (-1)^{|\lambda|} x_1^{-\lambda_m+m-1}\dots x_m^{-\lambda_1} y_1^{-\lambda'_n+n}\dots y_n^{-\lambda'_1+1},
 $$
one can check that 
$$
\sum_{w \in W_0} w \left( \frac{x_1^{-\lambda_m+m-1}\dots x_m^{-\lambda_1} y_1^{-\lambda'_n+n}\dots y_n^{-\lambda'_1+1}    }{\Delta(u)\Delta(v)\prod_{j=1}^n (y_j^{}-y_j^{-1})} \right)=0
$$
unless  $\lambda$ is either empty or partition of the form $$\lambda=(a_{1},\dots,a_{r}|a_{1}+1,\dots,a_{r}+1)$$ in Frobenius notations (see e.g. \cite{Mac}), in which case it is equal to $(-1)^{|\lambda|}.$
Now the proposition follows from the fact that the number of  such  diagrams contained in the $(m\times n)$ rectangle is equal to $2^{\min(m-1,n)},$ which can be proved by induction or reduced to the odd case.
\end{proof}

Let  $M=V^{\lambda}$ be the polynomial representation of $\mathfrak{gl}(m,n)$ determined by a Young diagram $\lambda \in H(m,n)$ and define the Euler supercharacters as
$
E_{\lambda} = E(V^{\lambda}).
$
Introduce the following polynomials 
\begin{equation}
\label{jac}
\varphi_a(z) = w^{a} + w^{-a} ,\: a>0,\:\varphi_{0}=1, \quad 
\psi_a(z) = \frac{w^{a+1} - w^{-a-1}}{w^{} - w^{-1}},
\end{equation}
where $z=w+w^{-1}.$ They are particular case of the Jacobi polynomials known as Chebyshev polynomials of the first and second kind respectively.

Let $i(\lambda), j(\lambda), k_i, l_j$ be defined by as before by (\ref{ilambda}),(\ref{jlambda}),(\ref{frob11}),(\ref{frob21}).
Define the following modification of $i(\lambda)$:
$$
i^{*}(\lambda)=\max\{i\mid\lambda_{i}+d-i> 0,\:\; 1\le i\le m\}.
$$
It is clear that that $i^*(\lambda)$ is equal to $i(\lambda)-1$ or to $i(\lambda)$ depending whether $l_{i(\lambda)}=0$ or not.

Similarly to Theorem 3.2  one can prove that

\begin{thm} If $l_{m}=0$ then the Euler supercharacters for $\mathfrak{osp}(2m,2n)$ can be given by the following formula
\begin{equation}\label{Euler11}
E_{\lambda}= C(\lambda)\sum_{w\in S_{m}\times S_{n}}w\left[\Pi_{\lambda}(u,v)\frac{\varphi_{l_{1}}(u_1)\dots \varphi_{l_{m}}(u_m)\psi_{k_{1}}(v_1)\dots \psi_{k_{n}}(v_n)}{\Delta(u)\Delta(v)}\right],
\end{equation}
where $C(\lambda) = (-1)^{b} 2^{m-i^*(\lambda)-1},\,\, b=\sum_{j>m}\lambda_{j}.$ 

If $l_{m}>0$ then we have a similar formula for the sum of Euler supercharacters
\begin{equation}\label{Euler21}
E_{\lambda}+ \theta(E_{\lambda})= (-1)^b\sum_{w\in S_{m}\times S_{n}}w\left[\Pi_{\lambda}(u,v)\frac{\varphi_{l_{1}}(u_1)\dots  \varphi_{l_{m}}(u_{m})\psi_{k_{1}}(v_1)\dots \psi_{k_{n}}(v_n)}{\Delta(u)\Delta(v)}\right].
\end{equation}
\end{thm}

As a corollary we have 
\begin{thm} \label{ej1}
The limit of the super Jacobi polynomials $SJ_{\lambda}(u,v,-1,p,q)$ as $(p,q) \rightarrow (0,0)$ is well defined and coincides up to a constant factor with the Euler supercharacter or sum of two Euler supercharacters:
\begin{equation}
\label{ej1}
\mathcal SJ_{\lambda}(u,v; -1,0,0)  =2^{i^*(\lambda)-m+1}E_{\lambda}(u,v)
\end{equation}
if $l_{m}=0$, and
\begin{equation}
\label{ej2}
\mathcal SJ_{\lambda}(u,v; -1,0,0) = E_{\lambda}+\theta(E_{\lambda})
\end{equation}
if $l_{m}>0.$ 
\end{thm}
The proof follows from comparison of formulas (\ref{Euler11}),(\ref{Euler21}) with (\ref{WeylJ}).

The Pieri formula for the Euler supercharacters of $\mathfrak{osp}(2m,2n)$  follows from Pieri formula for super Jacobi polynomials with $p=q=0$. 
Introduce the following functions:
\begin{equation}
\label{ad}
a_d(\mu,\lambda)=\begin{cases}0,\:\: \mu=\lambda\setminus\Box\,\, \text {and} \,\,j-i=-d  \cr 2,\,\,\, \mu=\lambda\setminus\Box\,\, \text{and}\,\, j-i=1-d \cr 1,\: \:\,\, \text{otherwise},\end{cases}
\end{equation}
As before $\mu\sim\lambda $ means that the Young diagram
$\mu$ can obtained from $\lambda$ by removing or adding one box.

\begin{thm}\label{Pieri1} Let $\mathcal SJ_{\lambda}(u,v; -1,0,0)=SJ_{\lambda}(u,v)$ be  specialized super Jacobi 
polynomials,  then the following Pieri formula holds
\begin{equation}\label{pieri1}
\left(\sum_{i=1}^m u_{i}-\sum_{j=1}^n v_{j}\right)\mathcal SJ_{\lambda}(u,v)=
\sum_{\mu \sim \lambda, \, \mu\in \mathcal P_{m,n}}a_d(\mu,\lambda)\mathcal SJ_{\mu}(u,v)
\end{equation}
\end{thm} 

The Jacobi-Trudy formula in this case follows directly from (\ref{JTinf}).

\begin{proposition} The specialized super Jacobi polynomials 
satisfy the following Jacobi--Trudy formula
\begin{equation}
\label{JTinfsupev}
\mathcal SJ_{\lambda}(u,v;-1,0,0)=\left|\begin{array}{cccc}
   h_{\lambda_{1}}&h^{(1)}_{\lambda_{1}}& \ldots &h^{(l-1)}_{\lambda_{1}}\\
  h_{\lambda_{2}-1}&h^{(1)}_{\lambda_{2}-1}& \ldots &h^{(l-1)}_{\lambda_{2}-1}\\
\vdots&\vdots&\ddots&\vdots\\
  h_{\lambda_{l}-l+1}&h^{(1)}_{\lambda_{l}-l+1}& \ldots &h^{(l-1)}_{\lambda_{l}-l+1}\\
 \end{array}\right|,
 \end{equation}
where $\lambda \in H_{m,n},\,d=m-n$  and $h^{(r)}_i$ are defined recursively by (\ref{relahinf}) with
$$
a(x)=0,\quad b(x)=1+\delta(x-1)-\delta(x)
$$
and $h_i^{(0)}=h_i=\mathcal SJ_{\lambda}(u,v;-1,0,0)$ for $\lambda=(i)$ for positive $i$ and $h_i^{(0)}=h_i\equiv 0$ if $i<0.$
\end{proposition}

Consider now a different  choice of Borel subalgebra and suitable parabolic subalgebras, following Gruson and Serganova \cite{GS}.

We assume for convenience  that $m> n$. Choose the set of simple roots with maximal number of isotropic roots:
$$
B=\{\varepsilon_{1}-\varepsilon_{2},\dots,\varepsilon_{m-n}-\delta_{1},\delta_{1}-\varepsilon_{m-n+1},\varepsilon_{m-n+1}-\delta_2,\dots, \delta_{n}-\varepsilon_{m},\delta_{n}+\varepsilon_{m}\}.
$$
The corresponding set of even positive roots is
$$
R^+_{0}=\{\varepsilon_{i}\pm\varepsilon_{j},\,\,i <j,\,\,\,\,\delta_{p}\pm\delta_{q},\,p<q,\,2\delta_{p}\}
$$
with the half-sum
$$
\rho_{0}=\frac12\sum_{\alpha\in R^+_0}\alpha=(m-1)\varepsilon_{1}+(m-2)\varepsilon_{2}\dots+\varepsilon_{m-1}+n\delta_{1}+(n-1)\delta_{2}+\dots+\delta_{n}.
$$
The set of positive odd roots is
$$
R^+_{1}=\{\varepsilon_{i}\pm\delta_{j},\,\,i -j< m-n,\,\,\,\,\delta_{j}\pm\varepsilon_{i},\,i-j\ge m-n\}
$$
with the half-sum
$$
\rho_{1}=n(\varepsilon_{1}+\dots+\varepsilon_{m-n})+(n-1)\varepsilon_{m-n+1}+\dots +\varepsilon_{m-1}+n\delta_{1}+\dots+\delta_{n},
$$
so
$$
\rho=\rho_{0}-\rho_{1}=\sum_{i=1}^{m-n}(m-n-i)\varepsilon_{i}
$$

\begin{lemma}\label{Young}  The weight 
$$
\chi=a_{1}\varepsilon_{1}+a_{2}\varepsilon_{2}+\dots+a_m\varepsilon_{m}+b_{1}\delta_{1}+b_{2}\delta_{2}+\dots+b_n\delta_{n}
$$
with $a_{m}\ge 0$ is a highest weight of an irreducible finite dimensional $\frak{osp}(2m,2n)$-module   if and only if there exists a partition $\lambda \in H(m,n)$  from the fat hook such that
\begin{equation}
\label{weight}
\chi+\rho= \sum_{i=1}^{i(\lambda)}(\lambda_{i}+d-i)\varepsilon_i +\sum_{j=1}^{j(\lambda)}(\lambda^{\prime}_{j}-d-j+1) \delta_j
\end{equation}
\end{lemma}
The proof again follows from comparison with Corollary 3 from \cite{GS}.
 
Let $\lambda$ be a partition from $H(m,n)$ and $\chi$ be the corresponding highest weight. Consider the parabolic subalgebra $\mathfrak p=\mathfrak p(\lambda)$ with
$$
R_{\mathfrak p}=\{\pm\varepsilon_{i}\pm\varepsilon_{j}, \,\,\pm\delta_{p}\pm\delta_{q},\,\,\pm2\delta_{p},\,  \pm \varepsilon_i \pm \delta_p\}
$$ 
with
$\,i(\lambda)\le i,j\le m,\,\,i\ne j,\,\,\,\,j(\lambda)\le p,q\le n,\,\,p\ne q\,,$ corresponding to the subalgebra $\mathfrak{osp}(2l,2l),\, l = m-i(\lambda)=n-j(\lambda).$

Let $E^{\mathfrak p}(\chi)$ be the corresponding Euler supercharacter given by (\ref{Vera}).
Let $\pi_{\lambda}$ be the same as in section 2 and define $t(\lambda)$ as the number of pairs $(i,j) \in \pi(\lambda)$ with $i-j\geq d=m-n.$ Similarly let $s(\lambda)$ be the number of pairs $(i,j) \in \lambda$ with $i-j\geq d.$

\begin{thm}
Let $\mathfrak p$ and $\chi$ be the parabolic subalgebra and its one-dimensional representation determined by a partition $\lambda \in H(m,n)$, then 
\begin{equation}
\label{ej11}
\mathcal SJ_{\lambda}(u,v; -1,0,0)  =(-1)^{s(\lambda)} 2^{i^*(\lambda)-i(\lambda)}E^{\mathfrak p}(\chi)
\end{equation}
if $\lambda_m \leq n$, and
\begin{equation}
\label{ej21}
\mathcal SJ_{\lambda}(u,v; -1,0,0) = (-1)^{s(\lambda)}(E^{\mathfrak p}(\chi) +\theta(E^{\mathfrak p}(\chi)))
\end{equation}
if $\lambda_m > n.$ 
\end{thm}

\begin{proof}
According to general formula (\ref{Vera})
$$E^{\mathfrak p}(\chi) = \sum_{w \in W_0} \, w\left(\frac{ \prod_{\alpha \notin R_{\mathfrak p} \cap R_1^+} (e^{\alpha/2} - e^{-\alpha/2})\, e^{\chi+\rho+\tau}}{\prod_{\alpha \in R_0^+}(e^{\alpha/2} - e^{-\alpha/2})}\right),
$$
where $$\tau = \frac{1}{2} \sum_{\alpha \in R_{\mathfrak p} \cap R_1^+} \alpha = \sum_{i>i(\lambda)}(m-i) \varepsilon_i + \sum_{j>j(\lambda)}(n-j+1) \delta_j.$$
One can check that
$$
 \prod_{\alpha \notin R_{\mathfrak p} \cap R_1^+} (e^{\alpha/2} - e^{-\alpha/2})= (-1)^{t(\lambda)}\prod_{(i,j)\in \pi_{\lambda}} (u_i-v_j)
$$
and that
$$
\chi+\rho+\tau= \sum_{i=1}^m l_i \varepsilon_i + \sum_{j=1}^{n} (k_j+1) \delta_j,$$
where $l_i,\, k_j$ are defined by (\ref{frob11}), (\ref{frob21}). Therefore we have
$$E^{\mathfrak p}(\chi)=(-1)^{t(\lambda)}\sum_{w \in W_0} \, w\left(\frac{\Pi_{\lambda}(u,v)e^{\chi+\rho+\tau}}{\Delta(u)\Delta(v)\prod_{j=1}^n (y_j-y_j^{-1})}\right)=
$$
$$(-1)^{t(\lambda)}\sum_{w \in W_0} \, w\left(\frac{\Pi_{\lambda}(u,v) \, x_1^{l_1}\dots x_m^{l_m} y_1^{k_1+1}\dots y_n^{k_n+1}}{\Delta(u)\Delta(v)}\right).
$$
If $l_m=0$ then it is easy to see that the average over $\mathbb Z_2^{m-1}\times \mathbb Z_2^{n}$ gives
$$E^{\mathfrak p}(\chi)= D(\lambda)\sum_{w \in S_m\times S_n} \, w\left(\frac{\Pi_{\lambda}(u,v) \, \varphi_{l_1}(u_1)\dots \varphi_{l_m}(u_m) \psi_{k_1+1}(v_1)\dots  \psi_{k_n+1}(v_n)}{\Delta(u)\Delta(v)}\right)
$$
with $D(\lambda)=(-1)^{t(\lambda)}2^{i(\lambda)-i^*(\lambda)}.$
If $l_m>0$ then we should consider the sum 
$$E^{\mathfrak p}(\chi) +\theta(E^{\mathfrak p}(\chi) = $$
$$
(-1)^{t(\lambda)}\sum_{w \in S_m\times S_n} \, w\left(\frac{\Pi_{\lambda}(u,v) \, \varphi_{l_1}(u_1)\dots \varphi_{l_m}(u_m) \psi_{k_1+1}(v_1)\dots  \psi_{k_n+1}(v_n)}{\Delta(u)\Delta(v)}\right).$$
Now the claim follows from Theorem \ref{main} and the relation
$$t(\lambda)=s(\lambda)+ b(\lambda).$$
\end{proof}

\begin{remark}
The coefficient $2^{i(\lambda)-i^*(\lambda)}$ (which can be either 1 or 2) can be eliminated by a different choice of parabolic subalgebra. Indeed, in the case when $i^*(\lambda)=i(\lambda)-1$ our choice of parabolic subalgebra $\mathfrak{p}$ is not the maximal one, which in that case corresponds to $\mathfrak{osp}(2l+2,2l).$
\end{remark}

\section{Concluding remarks}

We have shown that Euler supercharacters
for orthosymplectic Lie superalgebras coincide with specialized super Jacobi polynomials.
This fact seems to be very important and should have more conceptual proof.
One possibility is to show that the Euler supercharacters can be uniquely characterized by Pieri formula and the fact that they are the common eigenfunctions of the corresponding algebra of quantum integrals of deformed Calogero--Moser problem. From representation theory point of view this gives the action of the translation functors on Euler supercharacters. The results of Brundan \cite{Brun} show that one might be able to characterize in a similar way the irreducible characters.

Another interesting question is about  geometry of the genuine parameter space of the (super) Jacobi polynomials. We have seen that the specialization process requires some blowing-up procedure.
The calculations in the special case of $\mathfrak{osp}(3,2)$ indicate that this may lead to a description of the characters of important classes of finite-dimensional representations.

\section{Acknowledgements} We are  very grateful to Vera Serganova for useful and stimulating discussions of the results of \cite{GS} during the Algebraic Lie Theory programme at INI, Cambridge in March 2009.

This work has been partially supported by the EPSRC and by the
European Union through the FP6 Marie Curie RTN ENIGMA (contract
number MRTN-CT-2004-5652) and through ESF programme MISGAM.

\end{document}